\documentclass[11pt]{article}
\usepackage[a4paper,margin=27mm]{geometry}
\usepackage[T1]{fontenc}
\usepackage{lmodern}
\usepackage{amsmath,amssymb,amsthm,mathtools,mathrsfs}
\usepackage{microtype,booktabs,xcolor}
\usepackage[colorlinks=true,linkcolor=blue!55!black,citecolor=blue!55!black,urlcolor=blue!55!black]{hyperref}
\usepackage{aliascnt}
\usepackage[capitalize,noabbrev,nameinlink]{cleveref}
\numberwithin{equation}{section}
\allowdisplaybreaks[1]
\newtheorem{theorem}{Theorem}[section]
\newaliascnt{lemma}{theorem}
\newtheorem{lemma}[lemma]{Lemma}
\aliascntresetthe{lemma}
\newaliascnt{proposition}{theorem}
\newtheorem{proposition}[proposition]{Proposition}
\aliascntresetthe{proposition}
\newaliascnt{corollary}{theorem}
\newtheorem{corollary}[corollary]{Corollary}
\aliascntresetthe{corollary}
\theoremstyle{remark}
\newaliascnt{remark}{theorem}
\newtheorem{remark}[remark]{Remark}
\aliascntresetthe{remark}
\crefname{lemma}{lemma}{lemmas}
\crefname{proposition}{proposition}{propositions}
\crefname{corollary}{corollary}{corollaries}
\crefname{remark}{remark}{remarks}
\DeclareMathOperator{\spr}{spr}
\DeclareMathOperator{\Tr}{Tr}
\newcommand{\RR}{\mathbb R}
\newcommand{\CC}{\mathbb C}
\newcommand{\NN}{\mathbb N}
\newcommand{\Id}{\mathrm{Id}}
\newcommand{\OO}{\Omega^{(d)}_{\delta,\varepsilon}}
\newcommand{\bn}[2]{\binom{#1}{#2}}
\title{Strict Monotonicity of the First Positive Steklov Eigenvalue\\on Eccentric Spherical Shells in Every Dimension}
\author{Huaian Diao \qquad Jinrui Zhang\\[4pt]
\small School of Mathematics, Jilin University\\
\small Changchun 130024, P. R. of China\\[3pt]
\small \texttt{diao@jlu.edu.cn} \qquad \texttt{jerryzhang40@gmail.com}}
\date{}
\begin{document}
\maketitle
\begin{abstract}
For $d\ge2$, consider the region between two nested spheres of fixed radii, with the inner sphere translated a distance $\varepsilon$ from the outer center. We prove that the first positive eigenvalue of the Steklov problem, with the spectral boundary condition imposed on both spheres, decreases strictly with $\varepsilon$ throughout the noncontact interval. This settles the planar monotonicity conjecture formulated by Nigam and establishes its higher-dimensional counterparts. In dimension two we use a double-center holomorphic expansion and a pair of positive compact coefficient operators. For every $d\ge3$, solid Gegenbauer harmonics yield analogous axisymmetric and transverse operators, with a positive rank-one monopole correction. Weighted angular-mode comparison for $d\ge3$ and a direct matrix comparison for $d=2$ identify the first positive eigenvalue as the unique level-one crossing of an associated spectral radius. The first eigenspace has dimension $d$ in the concentric case and dimension one for every nonconcentric configuration. Convergence and reconstruction estimates justify the coefficient reductions at the two boundaries.
\end{abstract}
\noindent\textbf{Keywords.} Steklov eigenvalues; eccentric spherical shells; strict spectral monotonicity; Gegenbauer harmonics; positive compact operators; spectral radius.\\
\textbf{2020 Mathematics Subject Classification.} Primary 35P15; Secondary 58J50, 47B65, 35J25.

\section{Introduction}\label{sec:intro}
The Steklov eigenvalues of a smooth bounded domain are the eigenvalues
of its Dirichlet-to-Neumann operator and reflect the geometry of the
boundary through harmonic extension. They arise, for example, as spectral limits of models of vibrating
membranes whose mass concentrates near the boundary
\cite{GirouardPolterovich,ColboisGirouardGordonSher,DallaRivaProvenzano2018}.
Their dependence on domain geometry is also central to spectral
inequalities and shape optimization
\cite{BogoselBucurGiacomini2017}. We consider here how the low-lying
Steklov spectrum changes when an interior spherical obstacle is
translated while the radii of both boundary components remain fixed.

For $d\ge2$, $0<\delta<1$, and $0\le\varepsilon<1-\delta$, let
\begin{equation}\label{eq:domain}
 \Omega_{\delta,\varepsilon}^{(d)}
 =B_1(0)\setminus\overline{B_\delta(\varepsilon e_d)}
 \subset\RR^d,
 \qquad e_d=(0,\ldots,0,1).
\end{equation}
By rotational invariance, the displacement may be taken along $e_d$.
On this eccentric spherical shell, we study
\begin{equation}\label{eq:steklov}
 \begin{cases}
  \Delta u=0&\text{in }\Omega_{\delta,\varepsilon}^{(d)},\\
  \partial_\nu u=\lambda u&\text{on }\partial\Omega_{\delta,\varepsilon}^{(d)}.
 \end{cases}
\end{equation}
The Steklov condition is imposed on both boundary components, with
$\nu$ the outward unit normal of the shell. Its eigenvalues, counted
with multiplicity, satisfy
\begin{equation}\label{eq:spectrum}
 0=\sigma_0^{(d)}(\delta,\varepsilon)
 <\sigma_1^{(d)}(\delta,\varepsilon)
 \le\sigma_2^{(d)}(\delta,\varepsilon)\le\cdots\to\infty.
\end{equation}
The volume and boundary areas are independent of $\varepsilon$;
translation changes the relative positions of the two boundary
components and hence the harmonic extension between them. Our main
question concerns the dependence of the first positive eigenvalue
$\sigma_1^{(d)}(\delta,\varepsilon)$ on $\varepsilon$.

For fixed radii, the concentric shell uniquely maximizes the first
positive Steklov eigenvalue \cite[Theorem~1.1]{Ftouhi2022}:
\begin{equation}\label{eq:known-maximizer}
 \sigma_1^{(d)}(\delta,\varepsilon)
 <\sigma_1^{(d)}(\delta,0)
 \qquad(0<\varepsilon<1-\delta).
\end{equation}
This does not compare two positive displacements. Such a strict
monotonicity property was conjectured for planar eccentric annuli in
\cite[Conjecture~6.1]{Nigam2026} on the basis of numerical evidence.
Related results for eccentric domains concern mixed
Steklov--Dirichlet boundary conditions, with the Dirichlet condition
on the inner boundary
\cite{HongLimSeo2022,GavitonePiscitelli2024,GavitonePiscitelliCorrection2024}.
Higher-dimensional Fourier--Gegenbauer expansions have been used to
study near-contact behavior and coefficient approximations for such
mixed problems \cite{HongLimSeo2026,HongLeeLim2025}.
Second-order perturbations of Steklov eigenvalues under boundary
deformations have also been investigated \cite{AllandViator2026}.
In contrast, the family \eqref{eq:domain} retains the Steklov condition
on both boundary components; the degree-zero harmonic enters the
coefficient reduction through a rank-one term.

Write $\mathscr E_\lambda^{(d)}(\delta,\varepsilon)$ for the real
Steklov eigenspace of \eqref{eq:steklov} at $\lambda$ (see
\eqref{eq:eigenspace-definition}). Our main result is the following.
\begin{theorem}\label{thm:main}
For every integer $d\ge2$, $0<\delta<1$, and
$0\le\varepsilon_1<\varepsilon_2<1-\delta$,
\begin{equation}\label{eq:main-monotonicity}
 \sigma_1^{(d)}(\delta,\varepsilon_2)
 <\sigma_1^{(d)}(\delta,\varepsilon_1).
\end{equation}
Moreover,
\begin{equation}\label{eq:main-eigenspace-dim}
 \dim_{\mathbb R}
 \mathscr E^{(d)}_{\sigma_1^{(d)}(\delta,\varepsilon)}
 (\delta,\varepsilon)
 =\begin{cases}
 d,&\varepsilon=0,\\
 1,&0<\varepsilon<1-\delta.
 \end{cases}
\end{equation}
\end{theorem}

Thus the concentric maximality in \eqref{eq:known-maximizer}
extends to a strict ordering of all admissible displacements, and
the $d$-dimensional first eigenspace of the concentric shell becomes
one-dimensional for every nonzero displacement. The local analysis
in \Cref{prop:local-concentric-splitting} further gives a quadratic
eigenvalue splitting at $\varepsilon=0$.

The proof uses two-center harmonic expansions to reduce the boundary
problem to positive compact coefficient operators. Convergence and
reconstruction estimates establish the correspondence between their
eigenvectors and Steklov eigenfunctions. Spectral-radius comparison,
combined with an angular-mode argument in $d\ge3$, identifies the
first positive eigenvalue and proves its strict decrease.

\Cref{sec:prelim} gives the spectral and operator-theoretic
preliminaries. The harmonic reduction and its spectral comparison for
$d\ge3$ are developed in
\Cref{sec:multi,sec:analytic,sec:spectral}; the planar construction
appears in \Cref{sec:plane}. The common spectral characterization,
concentric case, and local eigenvalue splitting are treated in
\Cref{sec:all-dimensional}.

\section{Preliminaries}\label{sec:prelim}
The geometric hypotheses in \eqref{eq:domain} imply that
$\Omega_{\delta,\varepsilon}^{(d)}$ is connected, bounded, and has
$C^\infty$ boundary consisting of two disjoint spheres. We record
first the normalizations, followed by the spectral and
operator-theoretic facts used below.

For a displacement $h\in\RR^d$ with $|h|<1-\delta$, write
$\Omega_{\delta,h}^{(d)}=B_1(0)\setminus\overline{B_\delta(h)}$.
Set $\varepsilon=|h|$ and choose $Q\in O(d)$ such that
$Qh=\varepsilon e_d$ (with $Q=I_d$ if $h=0$). If $y=Qx$ and
$v(y)=u(Q^{\mathsf T}y)$, then
\begin{equation}\label{eq:rigid-covariance}
 \begin{aligned}
 \nabla_yv(y)&=Q\nabla_xu(x),&
 \Delta_yv(y)&=\Delta_xu(x),\\
 \nu_y(y)&=Q\nu_x(x),&
 \partial_{\nu_y}v(y)&=\partial_{\nu_x}u(x).
 \end{aligned}
\end{equation}
Consequently, with $\sigma_k(D)$ denoting the $k$th Steklov
eigenvalue of a smooth bounded domain $D$,
\begin{equation}\label{eq:rigid-isospectral}
 \sigma_k\bigl(\Omega_{\delta,h}^{(d)}\bigr)
 =\sigma_k\bigl(\Omega_{\delta,\varepsilon}^{(d)}\bigr),
 \qquad \varepsilon=|h|,\quad k\ge0.
\end{equation}
For a smooth bounded domain $\Omega\subset\RR^d$ and $R>0$,
\begin{equation}\label{eq:scaling-law}
 \sigma_k(R\Omega)=R^{-1}\sigma_k(\Omega),\qquad k\ge0.
\end{equation}
Thus the outer radius in \eqref{eq:domain} may be normalized to one.

\subsection{Steklov spectrum and variational estimates}
The spectral assertions in \eqref{eq:spectrum} are established for
smooth bounded connected domains in arbitrary dimension by the
Dirichlet-to-Neumann spectral theory
\cite[Sec.~1]{ColboisGirouardGordonSher}; see also
\cite[Sec.~1.1]{Ftouhi2022}. In particular, eigenfunctions belong
to $C^\infty(\overline\OO)$, since the boundary is $C^\infty$.
Let $\Tr:H^1(\OO)\to H^{1/2}(\partial\OO)$ be the Sobolev trace.
For $\lambda\ge0$, the real weak eigenspace is
\begin{equation}\label{eq:eigenspace-definition}
\begin{aligned}
 \mathscr E_{\lambda}^{(d)}(\delta,\varepsilon)
 :=\bigl\{u\in H^1(\OO;\mathbb R):\quad
 &\int_{\OO}\nabla u\cdot\nabla\varphi\,dx
 =\lambda\int_{\partial\OO}(\Tr u)(\Tr\varphi)\,dS,\\
 &\text{for every }\varphi\in H^1(\OO;\mathbb R)\bigr\}.
\end{aligned}
\end{equation}
Smooth-boundary regularity identifies this space with the classical
eigenspace of \eqref{eq:steklov}. By the trace theorem and compact
boundary embedding \cite[Secs.~1--2]{LambertiProvenzano2020},
\begin{equation}\label{eq:trace-compact}
 H^1(\OO)\xrightarrow{\,\Tr\,}H^{1/2}(\partial\OO)
 \hookrightarrow\hookrightarrow L^2(\partial\OO).
\end{equation}
For the connected shell, the first positive eigenvalue is characterized
by the attained variational minimum
\cite[Sec.~1.1]{Ftouhi2022}
\begin{equation}\label{eq:variational}
 \sigma_1^{(d)}(\delta,\varepsilon)
 =\min_{\substack{u\in H^1(\OO),\ \int_{\partial\OO}\Tr u\,dS=0\\
                   \|\Tr u\|_{L^2(\partial\OO)}=1}}
 \int_{\OO}|\nabla u|^2\,dx>0.
\end{equation}

\begin{lemma}\label{lem:prelim-bound}
For $d\ge2$, $0<\delta<1$, and $0\le\varepsilon<1-\delta$,
\begin{equation}\label{eq:range}
 0<\sigma_1^{(d)}(\delta,\varepsilon)
 \le\frac{1-\delta^d}{1+\delta^{d+1}}<1.
\end{equation}
\end{lemma}
\begin{proof}
Use $u(x)=x_1$ in \eqref{eq:variational}. Since both centers lie on
the $e_d$ axis, translation of the inner sphere gives
\begin{equation}\label{eq:test-mean}
 \int_{\partial\OO}\Tr(x_1)\,dS
 =\int_{\partial B_1(0)}x_1\,dS
 +\int_{\partial B_\delta(\varepsilon e_d)}x_1\,dS=0.
\end{equation}
If $\omega_d=|B_1(0)|$, direct integration gives
\begin{equation}\label{eq:test}
 \int_{\OO}|\nabla x_1|^2\,dx=\omega_d(1-\delta^d),\qquad
 \int_{\partial\OO}x_1^2\,dS=\omega_d(1+\delta^{d+1}).
\end{equation}
The conclusion follows from \eqref{eq:variational}.
\end{proof}
All coefficient operators below are constructed for $0\le\lambda<1$.
At $\varepsilon=0$, powers $\varepsilon^0$ in the
coefficient polynomials have their usual value $1$.

\subsection{Positive compact operators}
Throughout, $\mathbb N=\{1,2,\ldots\}$.
The space $\ell^2=\ell^2(\mathbb N;\mathbb R)$ carries its
coordinatewise order and the usual inner product.
A positive operator is order preserving; it need not be self-adjoint.
For a bounded real operator $T$, let $\sigma(T)$ denote the spectrum of
its complexification, and let $\Id$ be the identity operator. Its
spectral radius is
\begin{equation}\label{eq:spr-definition}
 \spr(T):=\sup\{\,|\mu|:\mu\in\sigma(T)\,\}.
\end{equation}
We use Gelfand's formula \cite[Ch.~III, \S3.2]{Kato1995},
\begin{equation}\label{eq:gelfand-formula}
 \spr(T)=\lim_{q\to\infty}\|T^q\|^{1/q}
          =\inf_{q\ge1}\|T^q\|^{1/q}.
\end{equation}
For a bounded operator $K$ with nonnegative matrix entries,
\begin{equation}\label{eq:diagonal-spr-bound}
 (K^q)_{11}\ge K_{11}^q,\qquad
 \spr(K)\ge\limsup_{q\to\infty}|(K^q)_{11}|^{1/q}
 \ge K_{11}.
\end{equation}

\begin{lemma}\label{lem:prelim-perron}
Let $K:\ell^2\to\ell^2$ be a positive compact operator with
$K_{nk}>0$ for every $n,k\ge1$. Then $\spr(K)>0$ is an
algebraically simple eigenvalue of $K$ and $K^*$.
Its right eigenvector can be chosen strictly positive in every
coordinate, and its left eigenvector can be chosen as
$w\in\ell^2$ with $w_n>0$ for every $n$.
\end{lemma}
\begin{proof}
We apply the Perron--Jentzsch theorem for ideal-irreducible positive
compact operators \cite[Theorem~1.5]{ChangWangWu2020}.
If $\mathcal J$ is a nonzero closed $K$-invariant order ideal and
$0\ne x\in\mathcal J$, then $|x|\in\mathcal J$ and
$K|x|\in\mathcal J$. As $K_{nk}>0$, every coordinate of $K|x|$
is positive, so
\begin{equation}\label{eq:irreducibility-check}
 0\le e_n\le\frac{K|x|}{(K|x|)_n},\qquad n\ge1.
\end{equation}
The ideal property gives $e_n\in\mathcal J$ for every $n$.
Density of finite sequences yields $\mathcal J=\ell^2$.
The cited theorem now gives the asserted eigenvalues and eigenvectors.
\end{proof}

\begin{lemma}\label{lem:prelim-radius-order}
Let $S,T$ be bounded operators on $\ell^2$ with nonnegative matrix
entries and $S_{nk}\le T_{nk}$ for all $n,k$. Suppose
$Sv=r_Sv$ for some $0\ne v\ge0$ and
$T^*w=r_Tw$ for some $w\in\ell^2$ with $w_n>0$ for every $n$.
If $(T-S)_{nk}>0$ for at least one pair $(n,k)$ with $v_k>0$,
then $r_T>r_S$.
In particular, when $r_S=\spr(S)$ and $r_T=\spr(T)$,
the spectral radii are strictly ordered.
\end{lemma}
\begin{proof}
For any bounded nonnegative matrix operator $U$ and
$x,w\in\ell^2$, monotone convergence and Cauchy--Schwarz give
\[
 \sum_{n,k\ge1}|w_n|U_{nk}|x_k|
 =\langle |w|,U|x|\rangle
 \le\|w\|_{\ell^2}\,\|U\|\,\|x\|_{\ell^2}<\infty.
\]
Thus the pairings below are absolutely convergent. Moreover,
$\langle w,v\rangle>0$, and
\[
 (r_T-r_S)\langle w,v\rangle
 =\langle w,(T-S)v\rangle>0.
\]
The strict inequality follows from the specified matrix entry.
\end{proof}

\begin{lemma}
\label{lem:prelim-radius-continuity}
Let $t\mapsto T_t$ be an operator-norm-continuous family of compact
operators on a complex Banach space. If $\spr(T_{t_0})>0$, then
\begin{equation}\label{eq:compact-radius-continuity}
 \lim_{t\to t_0}\spr(T_t)=\spr(T_{t_0}).
\end{equation}
\end{lemma}
\begin{proof}
Set $r_0=\spr(T_{t_0})>0$. Compactness gives an isolated
eigenvalue $\mu_0$ of finite algebraic multiplicity with
$|\mu_0|=r_0$. Under small norm perturbations, its Riesz
projection has locally constant positive rank, yielding
$\liminf_{t\to t_0}\spr(T_t)\ge r_0$
\cite[Ch.~IV, \S3, Thms.~3.15--3.16]{Kato1995}.
Upper semicontinuity of the spectrum gives
$\limsup_{t\to t_0}\spr(T_t)\le r_0$
\cite[Ch.~IV, \S3, Thm.~3.1]{Kato1995}.
Together these bounds imply \eqref{eq:compact-radius-continuity}.
\end{proof}

\section{Two-center harmonic expansions for \texorpdfstring{$d\ge3$}{d >= 3}}\label{sec:multi}
This section establishes the translation identities for the regular
and singular solid harmonics, then converts the two Steklov boundary
conditions \eqref{eq:steklov} into the coefficient matrices
\eqref{eq:Kd0}--\eqref{eq:Kd1}. The analytic equivalence of those
matrices with \eqref{eq:steklov} is proved separately in
\Cref{sec:analytic}.

Fix an integer $d\ge3$ and write
\begin{equation}\label{eq:dim-param}
 c=d-2\ge1,\qquad a=\frac c2,\qquad
 x=(x',\xi)\in\RR^{d-1}\times\RR,\quad r=|x|.
\end{equation}
The rotations of $x'$ preserve the geometry \eqref{eq:domain}.
Write $x=(r\sin\theta\,\omega,r\cos\theta)$, where
$\omega\in\mathbb S^{d-2}$. The angular Laplacian has the form
\begin{equation}\label{eq:angular-laplacian}
 \Delta_{\mathbb S^{d-1}}=\partial_{\theta}^2
 +(d-2)\cot\theta\,\partial_{\theta}
 +\frac1{\sin^2\theta}\Delta_{\mathbb S^{d-2}}.
\end{equation}
For a real-valued spherical harmonic $Y_m$ of degree $m$ on
$\mathbb S^{d-2}$,
\begin{equation}\label{eq:transverse-eigenvalue}
 -\Delta_{\mathbb S^{d-2}}Y_m=m(m+d-3)Y_m.
\end{equation}
Separation in \eqref{eq:angular-laplacian} yields the degree-$n$ angular harmonics
\begin{equation}\label{eq:full-angular-basis}
 \Psi_{n,m,Y}(\theta,\omega)
 = (\sin\theta)^m C_{n-m}^{a+m}(\cos\theta)Y_m(\omega),
 \qquad n\ge m\ge0,
\end{equation}
up to nonzero normalization constants (see \cite[Ch.~6, pp.~73--109]{AxlerBourdonRamey}).
Let $\Pi_m$ denote the $L^2(\mathbb S^{d-2})$-orthogonal projection
onto the entire degree-$m$ transverse spherical-harmonic space,
acting on the angular variable $\omega$. Gegenbauer orthogonality
and completeness give the sector decomposition
\eqref{eq:full-angular-basis}. In cylindrical coordinates
$x=(s\omega,\xi)$, $s>0$, the differential expression and the
Steklov normal derivative satisfy the following identities, where
$[A,B]:=AB-BA$ denotes the commutator:
\begin{equation}\label{eq:sector-projection-commutators}
 \begin{gathered}
 \Delta=\partial_s^2+\frac{d-2}{s}\partial_s
       +\partial_\xi^2+\frac1{s^2}\Delta_{\mathbb S^{d-2}},\\
 [\Pi_m,\Delta]u=0\quad\text{in }\OO,\qquad
 [\Pi_m,\partial_\nu]u=0\quad\text{on }\partial\OO,
 \end{gathered}
\end{equation}
for smooth $u$ (first away from $s=0$, then by continuity).
Smooth extension through the axis can be seen by expanding a smooth
function in Taylor polynomials in $x'$. Projection onto transverse
degree $m$ retains only terms of the form
$s^{m+2j}Y_m(\omega)$, which are polynomials in $x'$;
Taylor remainders of arbitrarily high order give the smooth
extension. This argument applies also to a fixed component
within that degree. The identities in
\eqref{eq:sector-projection-commutators} hold on the common
smooth domains of the respective operators. The boundary identity
uses the fact that both spherical normals have coefficients
independent of $\omega$.
The same decomposition holds for the weak problem, not only for
smooth functions. To see this, let $P$ be either $\Pi_m$ or the
orthogonal projection onto one fixed transverse harmonic of degree
$m$. For smooth functions, angular orthogonality and integration
by parts on $\mathbb S^{d-2}$ give
\[
 \begin{aligned}
 \|Pu\|_{L^2(\OO)}&\le\|u\|_{L^2(\OO)},&
 \|\nabla Pu\|_{L^2(\OO)}&\le\|\nabla u\|_{L^2(\OO)},\\
 \int_{\OO}\nabla(Pu)\cdot\nabla v\,dx
 &=\int_{\OO}\nabla u\cdot\nabla(Pv)\,dx,\\
 \int_{\partial\OO}\Tr(Pu)\Tr v\,dS
 &=\int_{\partial\OO}\Tr u\Tr(Pv)\,dS.
 \end{aligned}
\]
The physical-boundary measure has an angularly independent
factor, so the corresponding angular projection is a contraction
also on $L^2(\partial\OO)$. Smooth density and trace continuity
extend these identities to $H^1$ and show that $P$ commutes with
the trace. Its range is therefore a closed $H^1$ subspace.
Testing the weak equation for $u$ with $P\varphi$ proves that
$Pu$ satisfies that equation for every $\varphi\in H^1$.
In particular, every nonzero angular component of an eigenfunction
is an eigenfunction with the same eigenvalue. Completeness on
$\mathbb S^{d-2}$ and Fubini's theorem imply that a nonzero
$L^2(\OO)$ function has a nonzero component in at least one
fixed harmonic sector.
Consequently the Steklov quadratic form is block diagonal in $m$.
The $m=0,1$ sectors are represented by
\eqref{eq:Kd0}--\eqref{eq:Kd1}, while \eqref{eq:higher-angular}
excludes $m\ge2$.

\subsection{Gegenbauer harmonics and translation formulas}\label{subsec:gegenbauer}
Let $C_n^a$ denote the Gegenbauer polynomial, with the normalization of \cite[Eq.~(18.12.4)]{DLMF}, defined by
\begin{equation}\label{eq:gegenbauer-gf}
 (1-2tz+t^2)^{-a}=\sum_{n=0}^\infty C_n^a(z)t^n,
 \qquad |t|<1,\quad |z|\le1.
\end{equation}
We use the following unnormalized solid harmonics, where $\partial_1=\partial_{x_1}$:
\begin{align}
 R_n^{(0)}(x)&=r^n C_n^a(\xi/r),&
 I_n^{(0)}(x)&=r^{-n-c}C_n^a(\xi/r),&&n\ge0,\label{eq:multi-bases0}\\
 R_n^{(1)}(x)&=\partial_1 R_{n+1}^{(0)}(x),&
 I_n^{(1)}(x)&=\partial_1 I_{n-1}^{(0)}(x),&&n\ge1.
 \label{eq:multi-bases1}
\end{align}
For each $n\ge1$, the two functions in \eqref{eq:multi-bases1} have the same angular factor, with radial degrees $n$ and $-(n+c)$:
\begin{equation}\label{eq:associated}
 \begin{aligned}
 R_n^{(1)}(x)&=-c\,x_1r^{n-1}C_{n-1}^{a+1}(\xi/r),\\
 I_n^{(1)}(x)&=-c\,x_1r^{-n-c-1}C_{n-1}^{a+1}(\xi/r).
 \end{aligned}
\end{equation}
The lowest-degree modes fix the normalization in
\eqref{eq:multi-bases0}--\eqref{eq:multi-bases1}:
\begin{equation}\label{eq:multi-low-degree-normalization}
 \begin{gathered}
 R_0^{(0)}=1,\qquad R_1^{(0)}=c\xi,\qquad
 I_0^{(0)}=r^{-c},\\
 R_1^{(1)}=-c x_1,\qquad
 I_1^{(1)}=-c x_1r^{-c-2}.
 \end{gathered}
\end{equation}
The normalization in \eqref{eq:associated} can be checked directly
before applying any translation formula. Differentiating the regular
generating function with respect to $x_1$ gives
\[
 \sum_{k\ge0}t^k\partial_1R_k^{(0)}(x)
 =-c x_1t^2(1-2t\xi+t^2r^2)^{-a-1}.
\]
Its $t^{n+1}$ coefficient is the first line of
\eqref{eq:associated}. Similarly, differentiating
$|x-te_d|^{-c}$ gives
\[
 \sum_{k\ge0}t^k\partial_1I_k^{(0)}(x)
 =-c x_1r^{-c-2}
   \sum_{k\ge0}C_k^{a+1}(\xi/r)(t/r)^k,
 \qquad |t|<r.
\]
The $t^{n-1}$ coefficient is the second line of
\eqref{eq:associated}; these identities also agree with
\cite[Eq.~(18.9.19)]{DLMF}. The completeness assertion for the two bases is the spherical-harmonic decomposition in \eqref{eq:full-angular-basis}; the precise two-center version is proved in \Cref{lem:multi-split}.

We next give the exact translations needed to compare the boundaries centered at $0$ and at $\varepsilon e_d$.
\begin{lemma}\label{lem:multi-shift}
Fix an integer $d\ge3$, set $c=d-2$ as in \eqref{eq:dim-param},
and let $\varepsilon\ge0$. For $m\in\{0,1\}$ and
$k\ge m$,
\begin{align}
 R_k^{(m)}(x)
 &=\sum_{j=m}^{k}\bn{k+c+m-1}{k-j}
   \varepsilon^{k-j}R_j^{(m)}(x-\varepsilon e_d),
 \label{eq:reg-shift}\\
 I_j^{(m)}(x-\varepsilon e_d)
 &=\sum_{n=j}^\infty\bn{n-m}{j-m}
   \varepsilon^{n-j}I_n^{(m)}(x),
 \qquad |x|>\varepsilon.
 \label{eq:sing-shift}
\end{align}
The singular series converges locally uniformly with all derivatives for $|x|>\varepsilon$.
\end{lemma}
\begin{proof}
For the regular axial family, the generating function
\eqref{eq:gegenbauer-gf} gives
\begin{equation}\label{eq:regular-gen-shift}
 \begin{aligned}
 \sum_{k\ge0}t^kR_k^{(0)}(x)
 &=(1-2t\xi+t^2r^2)^{-a}\\
 &=(1-\varepsilon t)^{-c}
 \sum_{j\ge0}\left(\frac{t}{1-\varepsilon t}\right)^j
 R_j^{(0)}(x-\varepsilon e_d).
 \end{aligned}
\end{equation}
For $k\ge j\ge0$, the coefficient of $t^{k-j}$ in
$(1-\varepsilon t)^{-c-j}$ is
\[
 [t^{k-j}](1-\varepsilon t)^{-c-j}
 =\binom{k+c-1}{k-j}\varepsilon^{k-j}.
\]
Here $[t^k]$ denotes extraction of the coefficient of $t^k$.
Taking $[t^k]$ in \eqref{eq:regular-gen-shift} proves
\eqref{eq:reg-shift} for $m=0$.
For the singular axial family, the Newton generating formula is
\begin{equation}\label{eq:newton-gf}
 |x-\varepsilon e_d|^{-c}
 =\sum_{n\ge0}\varepsilon^n I_n^{(0)}(x),
 \qquad |x|>\varepsilon.
\end{equation}
Differentiating $j$ times with respect to $\varepsilon$ yields
\begin{equation}\label{eq:sing-deriv-shift}
 \begin{aligned}
 \partial_\varepsilon^j|x-\varepsilon e_d|^{-c}
 &=j!I_j^{(0)}(x-\varepsilon e_d),\\
 \partial_\varepsilon^j\sum_{n\ge0}\varepsilon^n I_n^{(0)}(x)
 &=j!\sum_{n\ge j}\binom nj\varepsilon^{n-j}I_n^{(0)}(x).
 \end{aligned}
\end{equation}
Equating the two lines in \eqref{eq:sing-deriv-shift} proves
\eqref{eq:sing-shift} for $m=0$.
Finally, apply $\partial_1$ to \eqref{eq:reg-shift} with $k$ replaced
by $k+1$, and to \eqref{eq:sing-shift} with $j$ replaced by $j-1$:
\begin{equation}\label{eq:transverse-shifts-derivation}
 \begin{aligned}
 \partial_1R_{k+1}^{(0)}(x)
 &=\sum_{j=1}^{k+1}\binom{k+c}{k+1-j}
 \varepsilon^{k+1-j}\partial_1R_j^{(0)}(x-\varepsilon e_d),\\
 \partial_1I_{j-1}^{(0)}(x-\varepsilon e_d)
 &=\sum_{n\ge j}\binom{n-1}{j-1}
 \varepsilon^{n-j}\partial_1I_{n-1}^{(0)}(x).
 \end{aligned}
\end{equation}
By \eqref{eq:multi-bases1}, \eqref{eq:transverse-shifts-derivation}
becomes \eqref{eq:reg-shift}--\eqref{eq:sing-shift} for $m=1$.
For each compact $E\Subset\{|x|>\varepsilon\}$, choose
$\varepsilon<\rho<\inf_{x\in E}|x|$. For every fixed derivative
order $N$, the differentiated generating series has a majorant
of the form
\begin{equation}\label{eq:shift-normal-majorant}
 \sum_{n\ge0}(n+1)^{N+c+2}
 \left(\frac{\rho}{\inf_{x\in E}|x|}\right)^n<\infty.
\end{equation}
The locally uniform convergence in \eqref{eq:shift-normal-majorant}
justifies the differentiations and rearrangements in
\eqref{eq:sing-deriv-shift}--\eqref{eq:transverse-shifts-derivation}.
\end{proof}

Before using the boundary condition \eqref{eq:steklov}, we
establish the unique decomposition \eqref{eq:multi-split}
and the coefficient regularity $A\in\ell^2$, required to give
meaning to the matrix equation \eqref{eq:multi-fixed}.
\begin{lemma}\label{lem:multi-split}
Fix $d\ge3$, $0<\delta<1$, $0\le\varepsilon<1-\delta$,
and $\lambda\ge0$. If
$u\in\mathscr E_{\lambda}^{(d)}(\delta,\varepsilon)$,
then the orthogonal projection $u_m$ onto the axial sector ($m=0$)
or a fixed transverse degree-one spherical-harmonic component
($m=1$) admits the following unique splitting. For $m=1$ choose
transverse coordinates so that the fixed component is proportional
to $\omega_1$; the corresponding statement for arbitrary
$Y_1(\omega)=b\cdot\omega$ is given by
\eqref{eq:transverse-b-normalization}. In these coordinates,
\begin{equation}\label{eq:multi-split}
 u_m(x)=\sum_{n\ge m} A_nR_n^{(m)}(x)
 +\sum_{j\ge m}B_jI_j^{(m)}(x-\varepsilon e_d),\qquad m=0,1,
\end{equation}
where the first series is harmonic in $B_1(0)$, the second is harmonic outside $\overline{B_\delta(\varepsilon e_d)}$ and decays at infinity, and the regular coefficients satisfy $A\in\ell^2$. Each series has a smooth one-sided trace on its corresponding physical boundary.
\end{lemma}
\begin{proof}
Fix auxiliary radii
\begin{equation}\label{eq:green-auxiliary}
 \delta<r_i,\qquad \varepsilon+r_i<r_o<1,
 \qquad D_{i,o}=B_{r_o}(0)\setminus\overline{B_{r_i}(\varepsilon e_d)}.
\end{equation}
With
$\Gamma_d(x)=[c|\mathbb S^{d-1}|]^{-1}|x|^{-c}$, and
$\delta_0$ denoting the Dirac mass at the origin, one has
$-\Delta\Gamma_d=\delta_0$ in the distributional sense. Green's representation (cf. \cite[Chs.~10--11, pp.~191--221]{AxlerBourdonRamey}) on $D_{i,o}$ reads
\begin{equation}\label{eq:green-decomposition}
 u_m(x)=\int_{\partial D_{i,o}}
 \left[\Gamma_d(x-y)\partial_{\nu_y}u_m(y)
       -u_m(y)\partial_{\nu_y}\Gamma_d(x-y)\right]dS_y,
 \qquad x\in D_{i,o}.
\end{equation}
The sign in \eqref{eq:green-decomposition} uses
$-\Delta_y\Gamma_d(x-y)=\delta_x(y)$: the second Green identity
has left-hand side $-u_m(x)$, which is moved to the other side.
Here $\nu_y=y/r_o$ on the outer sphere and
$\nu_y=-(y-\varepsilon e_d)/r_i$ on the inner sphere. Thus both
integrals use the outward normal of the shell, including its minus
sign on the inner component.
Let $F_{i,o}$ and $G_{i,o}$ be, respectively, the integrals in
\eqref{eq:green-decomposition} over the outer and inner spheres.
Because their kernels have no singularities on the opposite sides of
the integration spheres,
\begin{equation}\label{eq:green-domains}
 \Delta F_{i,o}=0\text{ in }B_{r_o}(0),\qquad
 \Delta G_{i,o}=0\text{ for }|x-\varepsilon e_d|>r_i,\qquad
 G_{i,o}(x)=O(|x|^{-c})\ (|x|\to\infty).
\end{equation}
The exact splitting $u_m=F_{i,o}+G_{i,o}$ holds on $D_{i,o}$.
Take two auxiliary pairs of radii
$(r_{i,1},r_{o,1})$ and $(r_{i,2},r_{o,2})$ with overlapping shells,
so that $\varepsilon+\max(r_{i,1},r_{i,2})
<\min(r_{o,1},r_{o,2})$. The identity
$F_1+G_1=F_2+G_2=u_m$ on the overlap defines
\begin{equation}\label{eq:green-gluing}
 W(x)=\begin{cases}
 F_1(x)-F_2(x),&|x|<\min(r_{o,1},r_{o,2}),\\
 G_2(x)-G_1(x),&|x-\varepsilon e_d|>
 \max(r_{i,1},r_{i,2}).
 \end{cases}
\end{equation}
The overlap identity $F_1-F_2=G_2-G_1$ and the auxiliary
separation inequality in \eqref{eq:green-auxiliary} give
\begin{equation}\label{eq:green-cover}
 \begin{aligned}
 F_1-F_2&=G_2-G_1
 &&\text{on }B_{r_*}(0)\cap\{|x-\varepsilon e_d|>\rho_*\},\\
 \RR^d&=B_{r_*}(0)\cup\{|x-\varepsilon e_d|>\rho_*\},
 &&r_*:=\min_i r_{o,i}>\varepsilon+\rho_*,\quad
   \rho_*:=\max(r_{i,1},r_{i,2}).
 \end{aligned}
\end{equation}
To check the covering, if $|x|\ge r_*$ then
$|x-\varepsilon e_d|\ge|x|-\varepsilon\ge r_*-\varepsilon>\rho_*$.
The two formulas agree on their open intersection, so they glue
as harmonic functions without a boundary jump. Thus
\eqref{eq:green-gluing} defines a harmonic function on all of
$\RR^d$. Its decay at infinity makes it bounded outside a ball;
continuity makes it bounded on that ball. By
\eqref{eq:green-domains} and the harmonic Liouville theorem,
\begin{equation}\label{eq:green-liouville}
 \Delta W=0\quad\text{in }\mathbb R^d,
 \qquad \lim_{|x|\to\infty}W(x)=0
 \quad\Longrightarrow\quad W\equiv0.
\end{equation}
Thus $F_1=F_2$ and $G_1=G_2$ on their common domains.
Two admissible pairs that do not satisfy the displayed overlap
inequality can each be compared with one common larger auxiliary
shell: choose its inner radius below both given inner radii, and
its outer radius above both given outer radii. Its overlap with
each original pair satisfies that inequality. Hence the decomposition
is independent of all admissible auxiliary choices.
Choose nested auxiliary shells with $r_i\downarrow\delta$ and
$r_o\uparrow1$; the unique compatible restrictions from
\eqref{eq:green-gluing}--\eqref{eq:green-liouville} define
functions in $\mathcal H(U)$, the space of real harmonic functions
on an open set $U$, satisfying
\begin{equation}\label{eq:green-global}
 F\in\mathcal H(B_1(0)),\qquad
 G\in\mathcal H(\mathbb R^d\setminus
                \overline{B_\delta(\varepsilon e_d)}),\qquad
 G(x)=O(|x|^{-c}),\qquad u_m=F+G.
\end{equation}
The decay of $G$ in \eqref{eq:green-global} follows from the
far-field estimate on any fixed auxiliary inner sphere.
Uniqueness also holds among all splittings with the stated
domains and decay, not only those produced by the integrals.
Indeed, for two such splittings the differences of the regular
parts and the opposite differences of the decaying parts agree
on $\OO$. The sets $B_1(0)$ and
$\{|x-\varepsilon e_d|>\delta\}$ cover $\RR^d$, since
$\varepsilon+\delta<1$. They therefore define an entire harmonic
function tending to zero at infinity, which vanishes by the
same Liouville argument. The decay condition excludes transferring
a nonzero constant between the two parts.

Recall the full transverse projector $\Pi_m$ from
\eqref{eq:sector-projection-commutators}. For $R\in O(d-1)$, define
$\mathcal U_Ru(s\omega,\xi):=u(sR^{-1}\omega,\xi)$.
The invariance of the degree-$m$ harmonic eigenspace under $O(d-1)$
gives the additional commutation identity
\begin{equation}\label{eq:green-projector-commutation}
 \Pi_m\mathcal U_Ru=\mathcal U_R\Pi_m u,
 \qquad R\in O(d-1).
\end{equation}
For a \emph{fixed} nonzero degree-one harmonic
$Y_1(\omega)=b\cdot\omega$, $b\in\RR^{d-1}\setminus\{0\}$,
define its rank-one angular projector, for $s>0$, by
\begin{equation}\label{eq:green-rank-one-projection}
 (\Pi_{Y_1}f)(s\omega,\xi)
 :=\frac{Y_1(\omega)}{\|Y_1\|_{L^2(\mathbb S^{d-2})}^2}
 \int_{\mathbb S^{d-2}}f(s\eta,\xi)Y_1(\eta)\,dS_\eta.
\end{equation}
Unlike $\Pi_1$, the projector $\Pi_{Y_1}$ need not commute with
all $\mathcal U_R$. It does commute with the Laplacian and the
physical normal derivative. Indeed, for $s>0$ and smooth $f$,
\begin{equation}\label{eq:green-rank-one-commutation}
 \begin{aligned}
 \Delta_{\mathbb S^{d-2}}\Pi_{Y_1}f
 &=-(d-2)\Pi_{Y_1}f
 =\Pi_{Y_1}\Delta_{\mathbb S^{d-2}}f,\\
 [\Pi_{Y_1},\partial_s]f
 &=[\Pi_{Y_1},\partial_\xi]f=0,\\
 [\Pi_{Y_1},\Delta]f&=0,\qquad
 [\Pi_{Y_1},\partial_\nu]f=0
 \quad\text{on the two spherical boundaries}.
 \end{aligned}
\end{equation}
The last identity uses the expression
$\partial_\nu=\nu_s(s,\xi)\partial_s+
\nu_\xi(s,\xi)\partial_\xi$ on either sphere.
The Green kernels in \eqref{eq:green-decomposition}, including
their normal derivatives on the auxiliary spheres, depend on
$\omega,\eta\in\mathbb S^{d-2}$ only through $\omega\cdot\eta$.
For any integrable zonal kernel $\Phi$, rotational symmetry gives
\begin{equation}\label{eq:green-degree-one-zonal}
 \int_{\mathbb S^{d-2}}\Phi(\omega\cdot\eta)
              (b\cdot\eta)\,dS_\eta
 = (b\cdot\omega)
   \int_{\mathbb S^{d-2}}\Phi(\omega\cdot\eta)
              (\omega\cdot\eta)\,dS_\eta .
\end{equation}
Indeed, the vector integral
$\int_{\mathbb S^{d-2}}\Phi(\omega\cdot\eta)\eta\,dS_\eta$
is parallel to $\omega$, and taking its scalar product with
$\omega$ yields \eqref{eq:green-degree-one-zonal}.
For $m=0$ the corresponding angular integral is independent of
$\omega$. In particular, when $u_m$ belongs to the fixed
$Y_1$-component, applying the angular projection to each integral
in \eqref{eq:green-decomposition} and using
\eqref{eq:green-degree-one-zonal} gives
\begin{equation}\label{eq:green-fixed-component-preservation}
 \Pi_{Y_1}u_m=u_m,\qquad
 \Pi_{Y_1}F_{i,o}=F_{i,o},\qquad
 \Pi_{Y_1}G_{i,o}=G_{i,o}.
\end{equation}
The same equalities hold with $\Pi_0$ in the axial sector.
Thus both Green integrals preserve each fixed degree-one
transverse harmonic component.
To identify the fixed component with the scalar matrices below,
let $b\ne0$ and choose $R_b\in O(d-1)$ so that
$R_bb=|b|e_1$. The associated solid harmonics are
\begin{equation}\label{eq:transverse-b-normalization}
 \begin{aligned}
 R_{n,b}^{(1)}&:=\sum_{i=1}^{d-1}b_i\partial_{x_i}R_{n+1}^{(0)},
 &I_{n,b}^{(1)}&:=\sum_{i=1}^{d-1}b_i\partial_{x_i}I_{n-1}^{(0)},\\
 \mathcal U_{R_b}R_{n,b}^{(1)}&=|b|R_n^{(1)},
 &\mathcal U_{R_b}I_{n,b}^{(1)}&=|b|I_n^{(1)}.
 \end{aligned}
\end{equation}
Because translations along $e_d$ commute with
$\mathcal U_{R_b}$, the coefficient identities
\eqref{eq:reg-shift}--\eqref{eq:sing-shift} hold with exactly
the same scalar coefficients for every nonzero $b$.
Combining \eqref{eq:green-global} and
\eqref{eq:transverse-b-normalization}, the regular and decaying
solid-harmonic expansions give \eqref{eq:multi-split} after this
harmless transverse rotation.
The coefficients in \eqref{eq:multi-split} are recovered by
spherical projections. Define
$\psi_n^{(m)}:=R_n^{(m)}|_{\mathbb S^{d-1}}$ as in
\eqref{eq:multi-bases0}--\eqref{eq:multi-bases1}. Orthogonality of
distinct spherical harmonic degrees gives, for $n,j\ge m$,
\begin{equation}\label{eq:two-center-coefficient-extraction}
 \begin{aligned}
 A_n&=\frac{1}{r^n\|\psi_n^{(m)}\|_{L^2(\mathbb S^{d-1})}^2}
     \int_{\mathbb S^{d-1}}F(r\widehat\omega)\psi_n^{(m)}(\widehat\omega)\,dS_{\widehat\omega},
       &&0<r<1,\\
 B_j&=\frac{\rho^{j+c}}{\|\psi_j^{(m)}\|_{L^2(\mathbb S^{d-1})}^2}
     \int_{\mathbb S^{d-1}}
     G(\varepsilon e_d+\rho\widehat\omega)\psi_j^{(m)}(\widehat\omega)\,dS_{\widehat\omega},
       &&\rho>\delta.
 \end{aligned}
\end{equation}
Here $\widehat\omega\in\mathbb S^{d-1}$ is the full-sphere integration
variable (distinct from $\omega\in\mathbb S^{d-2}$ in
\eqref{eq:full-angular-basis}), $F$ and $G$ are the functions in
\eqref{eq:green-global}, and the angular factors are evaluated
about their respective centers. In particular,
\[
 R_n^{(m)}(r\widehat\omega)=r^n\psi_n^{(m)}(\widehat\omega),
 \qquad
 I_n^{(m)}(\rho\widehat\omega)
 =\rho^{-n-c}\psi_n^{(m)}(\widehat\omega).
\]
The second identity uses the same angular factor, including its
sign and scale, in \eqref{eq:associated}. Both integrals in
\eqref{eq:two-center-coefficient-extraction} use unit-sphere
measure $dS_{\widehat\omega}$, yielding the stated radial
powers and squared norms. Combining
\eqref{eq:green-liouville} with
\eqref{eq:two-center-coefficient-extraction} proves uniqueness of
the full two-center coefficient representation.

To establish $A\in\ell^2$, we use the boundary regularity of the
harmonic terms. The function $G$ is harmonic
through the outer sphere, so $F=u_m-G$ has a $C^\infty$ trace there.
The functions $\psi_n^{(m)}$ introduced in
\eqref{eq:two-center-coefficient-extraction} are spherical harmonics
of degree $n$:
\begin{equation}\label{eq:angular-projection}
 -\Delta_{\mathbb S^{d-1}}\psi_n^{(m)}
  =n(n+c)\psi_n^{(m)},\qquad
 A_n=\frac{\langle F,\psi_n^{(m)}\rangle_{L^2(\mathbb S^{d-1})}}
 {\|\psi_n^{(m)}\|^2_{L^2(\mathbb S^{d-1})}}.
\end{equation}
The Gegenbauer norm identity \cite[Table~18.3.1, Gegenbauer row]{DLMF}
\begin{equation}\label{eq:gegenbauer-norm}
 \int_{-1}^1(1-t^2)^{\gamma-1/2}\bigl(C_n^{\gamma}(t)\bigr)^2dt
 =\frac{\pi2^{1-2\gamma}\Gamma(n+2\gamma)}
 {n!(n+\gamma)\Gamma(\gamma)^2},\qquad \gamma>0,
\end{equation}
For $m=0$, use $\gamma=a$ in \eqref{eq:gegenbauer-norm}; for $m=1$,
use $\gamma=a+1$ together with the factor $x_1$ in
\eqref{eq:associated}. The gamma-ratio bounds imply that, for some
$C_d\ge1$, $M_d\ge0$ independent of $n$,
\begin{equation}\label{eq:angular-polynomial-norm}
 C_d^{-1}(1+n)^{-M_d}
 \le \|\psi_n^{(m)}\|_{L^2(\mathbb S^{d-1})}
 \le C_d(1+n)^{M_d},\qquad m=0,1.
\end{equation}
For every $N\in\mathbb N$, integration by parts in
\eqref{eq:angular-projection} gives
\begin{equation}\label{eq:green-coefficient-decay}
 |A_n|\le\frac{\|(-\Delta_{\mathbb S^{d-1}})^N(F|_{\mathbb S^{d-1}})\|_{L^2}}
 {[n(n+c)]^N\|\psi_n^{(m)}\|_{L^2}}
 \le C_{N_0}(1+n)^{-N_0},\qquad n\ge\max(1,m)
\end{equation}
for any prescribed $N_0$ after choosing $N$ sufficiently large.
Thus $A\in\ell^2$. Finally $F$ is harmonic across the inner sphere,
so $G=u_m-F$ has a smooth one-sided trace there. Applying the
same spherical integration-by-parts argument to
$G(\varepsilon e_d+\delta\widehat\omega)$ shows that
$B_j\delta^{-j-c}$ decreases faster than every inverse power
of $j$. Angular differentiation of the basis and any fixed
number of radial derivatives introduce only polynomial factors
in the degree. Consequently the regular series and its
derivatives converge uniformly up to the outer sphere, and the
singular series and its derivatives converge uniformly from the
exterior up to the inner sphere. Projection at intermediate
radii and these one-sided limits justify the boundary coefficient
identities without assuming extension of $u_m$ through a
physical boundary.
\end{proof}

\subsection{Axial and transverse coefficient operators}\label{subsec:multi-operators}
The radial factors of the two multipole bases are $r^j$ and $r^{-j-c}$. This produces the following dimension-dependent Robin multipliers:
\begin{equation}\label{eq:multi-factors}
 \alpha_n^{(d)}(\lambda)=\frac{n+c+\lambda}{n-\lambda},\qquad
 \beta_j^{(d)}(\lambda)=\frac{j+\lambda\delta}{j+c-\lambda\delta},\qquad
 D_\lambda^{(d)}=c(1+\delta^{c+1})+\lambda(\delta^{c+1}-\delta).
\end{equation}
For $0\le\lambda<1$, $n,j\ge1$, and $c=d-2\ge1$, the factors
in \eqref{eq:multi-factors} are positive. In particular,
\begin{equation}\label{eq:multi-D-lower}
 D_\lambda^{(d)}
 \ge c(1+\delta^{c+1})-\delta
 >c-\delta>0.
\end{equation}

\begin{proposition}\label{prop:multi-operators}
Fix $d\ge3$, $0<\delta<1$, $0\le\varepsilon<1-\delta$,
and $0<\lambda<1$. For either $m\in\{0,1\}$, the
Steklov conditions \eqref{eq:steklov} imply the coefficient system
\begin{equation}\label{eq:multi-fixed}
 A=K^{(d,m)}_{\lambda,\varepsilon}A,
\end{equation}
with matrices indexed by $n,k\ge1$:
\begin{align}
 (K^{(d,0)}_{\lambda,\varepsilon})_{nk}
 &=\alpha_n^{(d)}\left[
 \sum_{j=1}^{\min(n,k)}\bn nj\bn{k+c-1}{k-j}
 \delta^{2j+c}\beta_j^{(d)}\varepsilon^{n+k-2j}
 +\frac{\lambda\delta^{c+1}}{D_\lambda^{(d)}}
 \bn{k+c-1}{k}\varepsilon^{n+k}\right],
 \label{eq:Kd0}\\
 (K^{(d,1)}_{\lambda,\varepsilon})_{nk}
 &=\alpha_n^{(d)}
 \sum_{j=1}^{\min(n,k)}\bn{n-1}{j-1}\bn{k+c}{k-j}
 \delta^{2j+c}\beta_j^{(d)}\varepsilon^{n+k-2j}.
 \label{eq:Kd1}
\end{align}
The last term of \eqref{eq:Kd0} defines a nonnegative operator
of rank one when $\lambda\varepsilon>0$, and the zero operator
otherwise. At $\lambda=0$, the formulas \eqref{eq:Kd0}--\eqref{eq:Kd1}
remain meaningful as coefficient operators, but the constant
Neumann mode is not encoded by $A=(A_n)_{n\ge1}$. The converse
for $0<\lambda<1$ is proved in \Cref{prop:multi-equivalence}.
\end{proposition}
\begin{proof}
Write the regular and singular series in \eqref{eq:multi-split} about the inner and outer centers, respectively. \Cref{lem:multi-shift} gives
\begin{equation}\label{eq:multi-HJ}
 H_j=\sum_{k\ge j}\bn{k+c+m-1}{k-j}\varepsilon^{k-j}A_k,
 \qquad
 J_n=\sum_{j=m}^{n}\bn{n-m}{j-m}\varepsilon^{n-j}B_j.
\end{equation}
For $j\ge1$, the inner condition $-\partial_\rho u=\lambda u$ and the outer condition $\partial_r u=\lambda u$ give, mode by mode,
\begin{equation}\label{eq:multi-BC}
 \begin{aligned}
 -jH_j\delta^{j-1}+(j+c)B_j\delta^{-j-c-1}
 &=\lambda(H_j\delta^j+B_j\delta^{-j-c}),\\
 nA_n-(n+c)J_n&=\lambda(A_n+J_n),
 \end{aligned}
\quad\Longleftrightarrow\quad
 \begin{cases}
 B_j=\delta^{2j+c}\beta_j^{(d)}H_j,\\
 A_n=\alpha_n^{(d)}J_n.
 \end{cases}
\end{equation}
For $m=0$, the $n=j=0$ projections give
\begin{equation}\label{eq:multi-zero}
 cB_0\delta^{-c-1}=\lambda(H_0+B_0\delta^{-c}),
 \qquad -cB_0=\lambda(A_0+B_0),
\end{equation}
where
\begin{equation}\label{eq:multi-H0}
 H_0=A_0+\sum_{k\ge1}\binom{k+c-1}{k}\varepsilon^kA_k,
 \qquad A_0=-\frac{c+\lambda}{\lambda}B_0\quad(\lambda>0).
\end{equation}
Substitute $A_0=-(c+\lambda)B_0/\lambda$ from
\eqref{eq:multi-H0} into the inner zero-mode equation in
\eqref{eq:multi-zero}. Before multiplying by $\delta^{c+1}$, this is
\begin{equation}\label{eq:zero-elimination-intermediate}
 \bigl[c\delta^{-c-1}+c+\lambda-\lambda\delta^{-c}\bigr]B_0
 =\lambda\sum_{k\ge1}\binom{k+c-1}{k}\varepsilon^kA_k.
\end{equation}
Multiplying \eqref{eq:zero-elimination-intermediate} by
$\delta^{c+1}$ gives
\begin{equation}\label{eq:zero-elimination}
 \left[c(1+\delta^{c+1})+\lambda(\delta^{c+1}-\delta)\right]B_0
 =\lambda\delta^{c+1}
 \sum_{k\ge1}\binom{k+c-1}{k}\varepsilon^kA_k.
\end{equation}
Since the bracket equals $D_\lambda^{(d)}>0$ by
\eqref{eq:multi-factors}, we obtain
\begin{equation}\label{eq:multi-B0}
 B_0=\frac{\lambda\delta^{c+1}}{D_\lambda^{(d)}}
 \sum_{k\ge1}\bn{k+c-1}{k}\varepsilon^kA_k,
 \qquad
 A_0=-\frac{(c+\lambda)\delta^{c+1}}{D_\lambda^{(d)}}
 \sum_{k\ge1}\bn{k+c-1}{k}\varepsilon^kA_k.
\end{equation}
The expression for $A_0$ in \eqref{eq:multi-B0} has cancelled the
factor $\lambda$ from $B_0$; for $\lambda>0$ it equals
$-(c+\lambda)B_0/\lambda$. This cancellation will distinguish a
continuous choice of constant from the free Neumann constant.
At $\lambda=0$, \eqref{eq:multi-zero} gives $B_0=0$, while $A_0$ is arbitrary. For $m=0$, substituting the two translations
\eqref{eq:multi-HJ} and the inner Robin coefficient in
\eqref{eq:multi-BC} into the outer equation gives
\begin{equation}\label{eq:multi-coefficient-sum}
 \begin{aligned}
 A_n=\alpha_n^{(d)}\Biggl[&\sum_{k\ge1}
 \biggl(\sum_{j=1}^{\min(n,k)}\binom nj
 \binom{k+c-1}{k-j}\delta^{2j+c}\beta_j^{(d)}
 \varepsilon^{n+k-2j}\biggr)A_k
 +\varepsilon^nB_0\Biggr].
 \end{aligned}
\end{equation}
Substitution of \eqref{eq:multi-B0} into
\eqref{eq:multi-coefficient-sum} gives \eqref{eq:Kd0}.
For $m=1$, the same substitution with the binomial factors
in \eqref{eq:multi-HJ} gives \eqref{eq:Kd1}.
For the forward implication, the sums in \eqref{eq:multi-HJ}
are first extracted at auxiliary interior spheres. For fixed $j$
and $0<\varepsilon<1$, the weights
$\binom{k+c+m-1}{k-j}\varepsilon^{k-j}$, $k\ge j$, belong
to $\ell^2$, since the binomial factor has polynomial growth
in $k$. Thus each $H_j$ series is absolutely convergent for
$A\in\ell^2$ by Cauchy--Schwarz, whereas each $J_n$ is a finite
sum. For $\varepsilon=0$, the translations are the identity.
One-sided smooth convergence from \Cref{lem:multi-split} permits
evaluation at the physical spheres. The sums over $j$ at each fixed
$n$ are finite; the double sums used in reconstruction are absolutely
convergent by \eqref{eq:multi-H-rebuild}--\eqref{eq:multi-J-rebuild}.
\end{proof}

\section{Compactness and reconstruction for \texorpdfstring{$d\ge3$}{d >= 3}}\label{sec:analytic}
The coefficient matrices \eqref{eq:Kd0}--\eqref{eq:Kd1}
are first derived for smooth solutions of \eqref{eq:steklov}.
\Cref{lem:multi-compact} proves that they act compactly on
$\ell^2$; \Cref{prop:multi-equivalence} then establishes the
reverse implication by normal convergence and boundary residual
identities \eqref{eq:inner-boundary-residual}--\eqref{eq:outer-boundary-residual}.

\subsection{Compactness and smoothing estimates}\label{subsec:multi-estimates}
We first dominate the positive matrix entries by a separable exponentially decaying kernel.
\begin{lemma}\label{lem:multi-compact}
Fix $d\ge3$, $0<\delta<1$, $0\le\lambda<1$ and
$0\le\varepsilon<1-\delta$. For $m\in\{0,1\}$, each
$K^{(d,m)}_{\lambda,\varepsilon}:\ell^2(\NN;\RR)
\to\ell^2(\NN;\RR)$ is a positive Hilbert--Schmidt operator. There exist constants $R>1$, $q\in(0,1)$, $C>0$, locally uniform in $(\lambda,\varepsilon)$, such that
\begin{equation}\label{eq:multi-HS-bound}
 0\le (K^{(d,m)}_{\lambda,\varepsilon})_{nk}
 \le C(k+1)^c R^{-n}q^k,
 \quad m=0,1,\quad n,k\ge1.
\end{equation}
Each fixed point $A=K^{(d,m)}_{\lambda,\varepsilon}A$ satisfies $|A_n|\le C_A R^{-n}$. Moreover, the operators are jointly Hilbert--Schmidt-norm continuous in $(\lambda,\varepsilon)$.
\end{lemma}
\begin{proof}
For $c=d-2\ge1$ and $1\le j\le k$,
\begin{equation}\label{eq:binom-majorant}
 \frac{\bn{k+c-1}{k-j}}{\bn kj}
 =\prod_{\ell=1}^{c-1}\frac{k+\ell}{j+\ell}
 \le C_c(k+1)^{c-1}.
\end{equation}
The empty product for $c=1$ equals one. Define
\begin{equation}\label{eq:multi-U}
 U^{(d)}_{nk}=\sum_{j=1}^{\min(n,k)}\bn nj\bn kj
 \delta^{2j+c}\varepsilon^{n+k-2j}.
\end{equation}
Define the analytic majorant map
\begin{equation}\label{eq:Tmap}
 T_\varepsilon(z)=\varepsilon+\frac{\delta^2z}{1-\varepsilon z},
 \qquad z\ge0,\quad \varepsilon z<1.
\end{equation}
For $j\ge1$ and $|t|<1$, the binomial series satisfies
\begin{equation}\label{eq:binom-series}
 \sum_{n\ge j}\binom nj t^{n-j}=(1-t)^{-j-1}.
\end{equation}
Apply \eqref{eq:binom-series} with $t=\varepsilon z$, multiply
by $\binom kj\delta^{2j+c}\varepsilon^{k-j}z^j$, and sum over
$1\le j\le k$. The binomial theorem gives
\begin{equation}\label{eq:multi-Ugen}
 \sum_{n\ge1}U^{(d)}_{nk}z^n
 =\frac{\delta^c}{1-\varepsilon z}
   \bigl(T_\varepsilon(z)^k-\varepsilon^k\bigr).
\end{equation}
The geometric restriction $\delta<1-\varepsilon$ gives
\begin{equation}\label{eq:multi-Rq-choice}
 T_\varepsilon(1)=\varepsilon+\frac{\delta^2}{1-\varepsilon}<1
 \quad\Longrightarrow\quad
 \exists R>1:\ \varepsilon R<1,\quad
 q:=T_\varepsilon(R)<1.
\end{equation}
Positivity of the coefficients in \eqref{eq:multi-Ugen} yields
\begin{equation}\label{eq:multi-Ubound}
 U^{(d)}_{nk}\le\frac{\delta^c}{1-\varepsilon R}R^{-n}q^k.
\end{equation}
The ratio of the $m=1$ and the $m=0$ binomial factors is
\begin{equation}\label{eq:multi-ratio}
 \frac{\bn{n-1}{j-1}\bn{k+c}{k-j}}
 {\bn nj\bn{k+c-1}{k-j}}
 =\frac{j(k+c)}{n(j+c)}\le k+c.
\end{equation}
The factors $\alpha_n^{(d)},\beta_j^{(d)}$ from \eqref{eq:multi-factors} are bounded on compact subsets of the parameter region. The coefficient $\bn{k+c-1}{k}$ in the rank-one term grows at most like $(k+1)^{c-1}$. Since $\varepsilon R<1$ and $\varepsilon\le q$, that term is bounded by $C(k+1)^{c-1}R^{-n}q^k$. Combining \eqref{eq:binom-majorant}--\eqref{eq:multi-ratio} gives \eqref{eq:multi-HS-bound}. The sum
\begin{equation}\label{eq:multi-HS-sum}
 \sum_{n,k\ge1}(k+1)^{2c}R^{-2n}q^{2k}<\infty
\end{equation}
proves the Hilbert--Schmidt property. For a fixed point, Cauchy--Schwarz gives
\begin{equation}\label{eq:multi-fixed-decay}
 |A_n|\le CR^{-n}\|A\|_{\ell^2}
 \left(\sum_{k\ge1}(k+1)^{2c}q^{2k}\right)^{1/2}.
\end{equation}
On any compact parameter set contained in
$\{0\le\lambda<1,\ 0\le\varepsilon<1-\delta\}$,
choose $R,q$ uniformly in \eqref{eq:multi-Ubound}.
The summable majorant \eqref{eq:multi-HS-sum} then permits
dominated convergence, proving continuity in Hilbert--Schmidt norm.
\end{proof}

The estimate \eqref{eq:multi-HS-bound} has a stronger consequence than
compactness: it maps arbitrary square-summable coefficient sequences
into an exponentially weighted sequence space. To formulate the
PDE--coefficient correspondence as a bounded map, we distinguish this
smoothing property from the regularity of an arbitrary harmonic
expansion.
\begin{proposition}\label{prop:multi-bounded-maps}
Fix $d\ge3$, $0<\delta<1$, $0\le\varepsilon<1-\delta$,
$m\in\{0,1\}$, and $0<\lambda<1$. For $m=1$ fix the transverse
direction of \eqref{eq:multi-bases1}; an arbitrary fixed nonzero
$Y_1$ is reduced to this direction by
\eqref{eq:transverse-b-normalization}. Define
\begin{equation}\label{eq:weighted-sequence-space}
 \mathscr X_R:=\left\{A=(A_n)_{n\ge1}:\
    \|A\|_{\mathscr X_R}:=\sup_{n\ge1}R^n|A_n|<\infty\right\},
 \qquad R>1.
\end{equation}
Write $\Pi_m^\sharp=\Pi_0$ for $m=0$ and
$\Pi_m^\sharp=\Pi_{Y_1}$ for $m=1$, with the projections of
\eqref{eq:sector-projection-commutators} and
\eqref{eq:green-rank-one-projection}. Let $\mathscr D'(\OO)$ denote the distributions on $\OO$, and define
\begin{equation}\label{eq:harmonic-H1-sector}
 \mathscr H_m^1(\OO):=\bigl\{u\in H^1(\OO;\RR):
 \Delta u=0\text{ in }\mathscr D'(\OO),\ 
 \Pi_m^\sharp u=u\bigr\}.
\end{equation}
Let $\mathcal C_m u=(A_n)_{n\ge1}$ be the regular coefficients
of the two-center Green decomposition
\eqref{eq:multi-split}. The following properties hold:
\begin{equation}\label{eq:bounded-extraction}
 \|\mathcal C_m u\|_{\ell^2}
 \le C_{d,\delta,\varepsilon,m}\|u\|_{H^1(\OO)},
 \qquad u\in\mathscr H_m^1(\OO).
\end{equation}
\begin{equation}\label{eq:bounded-regularization}
 K_{\lambda,\varepsilon}^{(d,m)}:\ell^2\longrightarrow\mathscr X_R,
 \qquad
 \|K_{\lambda,\varepsilon}^{(d,m)}A\|_{\mathscr X_R}
 \le C_{\lambda,\delta,\varepsilon,d}\|A\|_{\ell^2},
 \qquad A\in\ell^2.
\end{equation}
where $R>1$ is chosen as in \eqref{eq:multi-Rq-choice}.
Write $\mathbb N_0=\{0,1,2,\ldots\}$.
For this $R$, take $H_j$ from the first translation in
\eqref{eq:multi-HJ}, define the positive-degree $B_j$ using
only the inner Robin relation in \eqref{eq:multi-BC}, and,
when $m=0$, set $A_0,B_0$ by \eqref{eq:multi-B0}.
These prescriptions define a linear reconstruction map
\begin{equation}\label{eq:bounded-reconstruction}
 \mathcal S_{\lambda,m}:\mathscr X_R\longrightarrow
 C^\infty(\overline\OO)\cap\mathscr H_m^1(\OO),\qquad
 \|\mathcal S_{\lambda,m}A\|_{C^p(\overline\OO)}
 \le C_{p,\lambda,\delta,\varepsilon,d,R}\|A\|_{\mathscr X_R},
 \quad p\in\mathbb N_0.
\end{equation}
The maps satisfy
\begin{equation}\label{eq:bounded-left-inverse}
 \mathcal C_m\mathcal S_{\lambda,m}=\Id_{\mathscr X_R},
 \qquad
 \|\mathcal S_{\lambda,m}A\|_{H^1(\OO)}
 \le C\|A\|_{\ell^2}
 \quad\text{if }A=K_{\lambda,\varepsilon}^{(d,m)}A.
\end{equation}
The map $\mathcal S_{\lambda,m}$ enforces the inner boundary
condition for every $A\in\mathscr X_R$. The outer positive-degree
condition is $A_n=\alpha_n^{(d)}J_n$. For every $A\in\mathscr X_R$, the reconstructed
function $u=\mathcal S_{\lambda,m}A$ satisfies the entire inner
boundary condition. On the outer boundary its residual is
\[
 (\partial_\nu u-\lambda u)|_{|x|=1}
 =\sum_{n\ge1}(n-\lambda)
 \bigl(A_n-(K^{(d,m)}_{\lambda,\varepsilon}A)_n\bigr)
 \psi_n^{(m)}.
\]
This series converges with every angular derivative. Its axial
zero mode is zero by \eqref{eq:multi-B0}. Thus reconstruction
solves both boundary conditions precisely when $A$ is a fixed point.
The bound in \eqref{eq:bounded-reconstruction} is asserted on
$\mathscr X_R$, not on the whole of $\ell^2$. In particular,
\begin{equation}\label{eq:bounded-composed-map}
 \mathcal S_{\lambda,m}K_{\lambda,\varepsilon}^{(d,m)}:
 \ell^2\longrightarrow H^1(\OO)
 \quad\text{is bounded}.
\end{equation}
\end{proposition}
\begin{proof}
We first show \eqref{eq:bounded-extraction}, including the dependence
on the unnormalized solid-harmonic basis. Choose an auxiliary
sphere $\Gamma_i=\partial B_{r_i}(\varepsilon e_d)$ as in
\eqref{eq:green-auxiliary} and an open neighborhood
$V_i\Subset\OO$ of $\Gamma_i$.
Interior elliptic regularity for the harmonic function $u$ and
the ordinary trace estimate on $\Gamma_i$ yield
\begin{equation}\label{eq:green-interior-bound}
 \|u\|_{L^2(\Gamma_i)}+
 \|\partial_\nu u\|_{L^2(\Gamma_i)}
 \le C\|u\|_{H^2(V_i)}
 \le C\|u\|_{H^1(\OO)}.
\end{equation}
In the inner Green potential $G$ of
\eqref{eq:green-decomposition}, the kernels and all their
$x$-derivatives are bounded on
$\partial B_1(0)\times\Gamma_i$, since
$\operatorname{dist}(\partial B_1(0),\Gamma_i)>0$.
Thus \eqref{eq:green-interior-bound}, the outer trace bound, and
$F=u-G$ in a collar of $\partial B_1(0)$ give
\begin{equation}\label{eq:green-outer-Hhalf}
 \|G\|_{H^{1/2}(\partial B_1(0))}
 +\|F\|_{H^{1/2}(\partial B_1(0))}
 \le C\|u\|_{H^1(\OO)}.
\end{equation}
The traces in \eqref{eq:green-outer-Hhalf} are well defined:
$G$ is smooth in that collar and $u$ is $H^1$ there, so $F$ has
an $H^{1/2}$ trace. In the remaining compact subball $F$ is
harmonic and smooth; hence $F\in H^1(B_1(0))$. The trace is
therefore also the boundary limit of its regular harmonic expansion.
Here the $H^{1/2}$ norm on the unit sphere can be defined
through the eigenvalues $n(n+c)$ of
$-\Delta_{\mathbb S^{d-1}}$. In particular, for the angular
functions $\psi_n^{(m)}$ of
\eqref{eq:two-center-coefficient-extraction},
\begin{equation}\label{eq:angular-norm-exact}
 \begin{aligned}
 \|\psi_n^{(0)}\|_{L^2(\mathbb S^{d-1})}^2
 &=|\mathbb S^{d-2}|\,h_n(a),\\
 \|\psi_n^{(1)}\|_{L^2(\mathbb S^{d-1})}^2
 &=\frac{c^2|\mathbb S^{d-2}|}{d-1}\,h_{n-1}(a+1),\\
 h_j(\gamma)&:=\frac{\pi2^{1-2\gamma}\Gamma(j+2\gamma)}
 {j!(j+\gamma)\Gamma(\gamma)^2},\qquad j\ge0,\ \gamma>0.
 \end{aligned}
\end{equation}
Indeed, write
$\widehat\omega=(\sqrt{1-t^2}\,\omega,t)$ with
$\omega\in\mathbb S^{d-2}$ and $-1<t<1$. The sphere measure is
\[
 dS_{\widehat\omega}=(1-t^2)^{(d-3)/2}\,dt\,dS_\omega,
 \qquad
 \psi_n^{(1)}(\widehat\omega)
 =-c\sqrt{1-t^2}\,\omega_1 C_{n-1}^{a+1}(t).
\]
Thus the axial squared norm has weight $(1-t^2)^{a-1/2}$,
and the transverse squared norm has weight
$(1-t^2)^{a+1/2}$ together with the factor $c^2\omega_1^2$.
Integrating first in $\omega$, using
$\int_{\mathbb S^{d-2}}\omega_1^2\,dS
=|\mathbb S^{d-2}|/(d-1)$, and then applying
\eqref{eq:gegenbauer-norm} proves both exact norms. Write $X_n\asymp Y_n$ when there exist constants
$0<c_1\le c_2<\infty$, independent of $n$, such that
$c_1Y_n\le X_n\le c_2Y_n$. For example, the degree-zero and degree-one norms are exactly
\[
 \|\psi_0^{(0)}\|_{L^2}^2=|\mathbb S^{d-1}|,\qquad
 \|\psi_1^{(0)}\|_{L^2}^2
 =\|\psi_1^{(1)}\|_{L^2}^2
 =\frac{c^2|\mathbb S^{d-1}|}{d},
\]
consistent with \eqref{eq:multi-low-degree-normalization}.
Cancellation of the Gamma factors also gives the exact ratio
\[
 \frac{\|\psi_n^{(1)}\|_{L^2}^2}{\|\psi_n^{(0)}\|_{L^2}^2}
 =\frac{n(n+c)}{c+1},\qquad n\ge1.
\]
The Gamma-ratio estimates yield
\begin{equation}\label{eq:angular-norm-weights}
 \|\psi_n^{(0)}\|_{L^2}^2\asymp(n+1)^{c-2},\qquad
 \|\psi_n^{(1)}\|_{L^2}^2\asymp(n+1)^c,
 \qquad
 (1+n(n+c))^{1/2}\|\psi_n^{(m)}\|_{L^2}^2\ge c_{d,m}>0.
\end{equation}
In the borderline case $d=3$ (so $c=1$), these estimates can
be read without asymptotic notation:
\[
 \|\psi_n^{(0)}\|_{L^2(\mathbb S^2)}^2=\frac{4\pi}{2n+1},\qquad
 \|\psi_n^{(1)}\|_{L^2(\mathbb S^2)}^2
 =\frac{2\pi n(n+1)}{2n+1},\qquad n\ge1.
\]
In particular,
$\sqrt{1+n(n+1)}\,4\pi/(2n+1)\ge2\pi$.
Thus the axial extraction estimate at $d=3$ uses exactly the
available $H^{1/2}$ trace regularity and loses no positive power
of $n$.
By orthogonality across the spherical-harmonic degrees
and the spectral $H^{1/2}$ norm,
\begin{equation}\label{eq:extraction-energy-bound}
 \begin{aligned}
 \|\mathcal C_m u\|_{\ell^2}^2
 &\le c_{d,m}^{-1}\sum_{n\ge1}
 (1+n(n+c))^{1/2}|A_n|^2\|\psi_n^{(m)}\|_{L^2}^2\\
 &\le C\|F|_{\partial B_1(0)}\|_{H^{1/2}}^2
 \le C\|u\|_{H^1(\OO)}^2,
 \end{aligned}
\end{equation}
which is \eqref{eq:bounded-extraction}. The decomposition used
here exists for every $u\in\mathscr H_m^1(\OO)$:
the auxiliary spheres lie strictly inside $\OO$, so the Green
representation and the uniqueness/gluing argument
\eqref{eq:green-decomposition}--\eqref{eq:green-global}
require only interior smoothness of harmonic $H^1$ functions.

By \eqref{eq:multi-HS-bound} and Cauchy--Schwarz,
\begin{equation}\label{eq:regularization-Cauchy}
 \sup_{n\ge1}R^n|(K_{\lambda,\varepsilon}^{(d,m)}A)_n|
 \le C\sum_{k\ge1}(k+1)^cq^k|A_k|
 \le C\left(\sum_{k\ge1}(k+1)^{2c}q^{2k}\right)^{1/2}
 \|A\|_{\ell^2},
\end{equation}
proving \eqref{eq:bounded-regularization}.
For $A\in\mathscr X_R$ choose $1<S<R$; then
$|A_k|\le\|A\|_{\mathscr X_R}S^{-k}$. Using
\eqref{eq:multi-HJ} and the binomial series gives
\begin{equation}\label{eq:bounded-rebuild-HB}
 \begin{aligned}
 |H_j|&\le C\|A\|_{\mathscr X_R}(S-\varepsilon)^{-j},\\
 |B_j|&\le C\|A\|_{\mathscr X_R}
 \left(\frac{\delta^2}{S-\varepsilon}\right)^j,
 \qquad j\ge1.
 \end{aligned}
\end{equation}
The degree-zero coefficients are bounded by
\eqref{eq:multi-B0} (when $m=0$).
Since $S>1>\varepsilon+\delta$, we have
$\rho_0:=\delta^2/(S-\varepsilon)<\delta$.
The solid-harmonic derivative bounds
\eqref{eq:multi-harmonic-derivative-bound}, applied with
$\sup_{\overline\OO}|x|<S$ and
$\inf_{\overline\OO}|x-\varepsilon e_d|=\delta>\rho_0$,
therefore imply, for every $p\ge0$,
\begin{equation}\label{eq:bounded-rebuild-series}
 \|\mathcal S_{\lambda,m}A\|_{C^p(\overline\OO)}
 \le C_p\|A\|_{\mathscr X_R}
 \left[\sum_{n\ge1}(n+1)^{M_p}S^{-n}
 +\sum_{j\ge1}(j+1)^{M_p}
 \left(\frac{\rho_0}{\delta}\right)^j+1\right]
 <\infty.
\end{equation}
The sum in \eqref{eq:bounded-rebuild-series} is harmonic and has
its indicated $m$-sector; hence
\eqref{eq:bounded-reconstruction} follows. By uniqueness of
the Green decomposition in \eqref{eq:green-global}, its regular
coefficient sequence is the prescribed $A$, so
$\mathcal C_m\mathcal S_{\lambda,m}=\Id_{\mathscr X_R}$.
If $A=K^{(d,m)}_{\lambda,\varepsilon}A$, then
\eqref{eq:bounded-regularization} gives
$\|A\|_{\mathscr X_R}\le C\|A\|_{\ell^2}$;
combining this bound with
\eqref{eq:bounded-reconstruction} proves the final inequality
in \eqref{eq:bounded-left-inverse}. Applying the same two
bounds without the fixed-point assumption proves
\eqref{eq:bounded-composed-map}.
Finally, writing $K=K^{(d,m)}_{\lambda,\varepsilon}$,
substitution of \eqref{eq:multi-HJ} and
\eqref{eq:multi-BC} gives $\alpha_n^{(d)}J_n=(KA)_n$ for the
reconstructed coefficients, with the $B_0$ term supplied by
\eqref{eq:multi-B0}. The outer residual is consequently
$(n-\lambda)(A_n-(KA)_n)$ in degree $n\ge1$.
The inner residual and both axial zero residuals vanish by
construction. The normal-convergence estimates above justify
these statements as equalities of smooth boundary functions.
\end{proof}

\begin{remark}\label{rem:ell2-versus-H1}
The $\ell^2$ norm of regular coefficients alone does not control
the $H^1$ energy of the associated regular-harmonic series.
For the degree-$n$ transverse regular harmonic,
homogeneity and Green's identity give
\begin{equation}\label{eq:regular-solid-energy}
 \int_{B_1(0)}|\nabla R_n^{(1)}|^2\,dx
 =n\|\psi_n^{(1)}\|_{L^2(\mathbb S^{d-1})}^2
 \asymp n^{c+1}\xrightarrow[n\to\infty]{}\infty,
\end{equation}
whereas its regular coefficient sequence has $\ell^2$ norm one.
Since $B_\delta(\varepsilon e_d)\subset B_{\varepsilon+\delta}(0)$
and $\varepsilon+\delta<1$, the solid-harmonic derivative bounds yield
\begin{equation}\label{eq:regular-shell-energy}
 \begin{aligned}
 0&\le\int_{B_\delta(\varepsilon e_d)}|\nabla R_n^{(1)}|^2\,dx
 \le C(n+1)^M(\varepsilon+\delta)^{2n-2}\longrightarrow0,\\
 &\hspace{12mm}\int_{\OO}|\nabla R_n^{(1)}|^2\,dx\asymp n^{c+1}.
 \end{aligned}
\end{equation}
Thus the reconstruction estimate in
\eqref{eq:bounded-reconstruction} is formulated on $\mathscr X_R$;
its bounded $\ell^2$ composition is
$\mathcal S_{\lambda,m}K^{(d,m)}_{\lambda,\varepsilon}$,
see \eqref{eq:bounded-composed-map}.

For $m=0$, the zero-mode equations at $\lambda=0$ give $B_0=0$
and leave the constant coefficient $A_0$ free. The continuous
extension of the cancelled formula \eqref{eq:multi-B0} selects
\[
 A_0=-\frac{\delta^{c+1}}{1+\delta^{c+1}}
       \sum_{k\ge1}\binom{k+c-1}{k}\varepsilon^k A_k.
\]
Choosing $A_0=0$ instead, as in
\eqref{eq:multi-neumann-reconstruction}, changes the reconstructed
function only by a constant. Hence the positive-eigenvalue
isomorphism of \Cref{prop:multi-equivalence} does not extend to the
Neumann eigenspace. More generally, $\mathcal C_m$ need not be
injective on all of $\mathscr H_m^1(\OO)$, since a singular multipole
may have vanishing regular coefficients. For $m=1$ no constant
coefficient occurs.

The reconstruction estimates are locally uniform for fixed $d,\delta$
on rectangles $0\le\lambda\le\Lambda<1$ and
$0\le\varepsilon\le E<1-\delta$: choose common radii and
majorants with
\[
 \delta<r_i<1-E,\qquad E+r_i<r_o<1,\qquad
 R>1,\quad ER<1,\quad q=T_E(R)<1,\quad 1<S<R.
\]
These choices provide uniform separation of auxiliary spheres and
uniform geometric-series bounds, including at $\lambda=0$ by the
cancelled zero-mode formula. The estimates concern the $H^1$ norms
on each of the varying shells.
\end{remark}

\subsection{Reconstruction and spectral equivalence}\label{subsec:multi-equivalence}
The bounded maps in \Cref{prop:multi-bounded-maps}, combined
with the translation formulas \eqref{eq:multi-HJ}, yield the
normal-convergence bounds \eqref{eq:multi-H-rebuild} and
\eqref{eq:multi-B-rebuild}, which justify evaluating normal
derivatives at both boundaries of \eqref{eq:domain}.
\begin{proposition}\label{prop:multi-equivalence}
Fix $d\ge3$, $0<\delta<1$, $0\le\varepsilon<1-\delta$,
$m\in\{0,1\}$ and $0<\lambda<1$. Then
\begin{equation}\label{eq:multi-equiv}
 \lambda\text{ is a Steklov eigenvalue in sector }m
 \ \Longleftrightarrow\ 
 \ker(\Id-K^{(d,m)}_{\lambda,\varepsilon})\ne\{0\}.
\end{equation}
For $m=0$ this is an isomorphism of the axial Steklov
eigenspace with the displayed kernel. For $m=1$ it is an
isomorphism \emph{for each fixed nonzero transverse harmonic}
$Y_1(\omega)=b\cdot\omega$; it does not identify the entire
degree-one transverse eigenspace with one copy of $\ell^2$.
For $\lambda=0$, any nonzero fixed point of
$K^{(d,0)}_{0,\varepsilon}$ reconstructs a nonconstant harmonic
function with homogeneous Neumann conditions.
\end{proposition}
\begin{proof}
For $0<\lambda<1$, \Cref{lem:multi-split} and the coefficient
identities \eqref{eq:multi-BC}--\eqref{eq:multi-B0} yield the map
\[
 u_m\quad(\lambda\text{-eigenfunction in sector }m)
 \longmapsto A\in\ker(\Id-K^{(d,m)}_{\lambda,\varepsilon}).
\]
The regular coefficient vector is nonzero for every nonzero
$u_m$: by \eqref{eq:multi-HJ}--\eqref{eq:multi-B0},
\[
 A=0\Longrightarrow H_j=B_j=0\quad(j\ge1),
 \qquad A_0=B_0=0\quad(m=0).
\]

Conversely, suppose $A\in\ell^2$ is a nonzero fixed point. By \eqref{eq:multi-fixed-decay}, choose $S$ with $1<S<R$ and $S-\varepsilon>\delta$ such that $|A_k|\le C S^{-k}$. The regular translations \eqref{eq:multi-HJ} then satisfy, with $m=0$ or $1$,
\begin{equation}\label{eq:multi-H-rebuild}
 \begin{aligned}
 |H_j|
 &\le C\sum_{k\ge j}\bn{k+c+m-1}{k-j}
              \varepsilon^{k-j}S^{-k}\\
 &= C\left(\frac S{S-\varepsilon}\right)^{c+m}
          (S-\varepsilon)^{-j}.
 \end{aligned}
\end{equation}
Using the boundedness of $\beta_j^{(d)}(\lambda)$ for $j\ge1$ gives
\begin{equation}\label{eq:multi-B-rebuild}
 |B_j|\le C_{d,\delta,\lambda,S}\,
       \delta^c\left(\frac{\delta^2}{S-\varepsilon}\right)^j,
 \qquad \rho_0:=\frac{\delta^2}{S-\varepsilon}<\delta.
\end{equation}
For $m=0$, define $A_0,B_0$ by \eqref{eq:multi-B0}; no
degree-zero coefficients occur when $m=1$.
Set $\theta:=\varepsilon+\rho_0<\varepsilon+\delta<1$.
The second translation in \eqref{eq:multi-HJ}, together with
\eqref{eq:multi-B-rebuild}, gives the additional estimate
\begin{equation}\label{eq:multi-J-rebuild}
 \begin{aligned}
 |J_n|&\le \sum_{j=m}^n\binom{n-m}{j-m}
             \varepsilon^{n-j}|B_j|,\\
 |J_n|&\le C_{d,\delta,\lambda,S}
 \begin{cases}
 (\varepsilon+\rho_0)^n,&m=0,\\
 \rho_0(\varepsilon+\rho_0)^{n-1},&m=1,
 \end{cases}\\
 &\hspace{32mm}n\ge\max(1,m).
 \end{aligned}
\end{equation}
For $m=0$, the $j=0$ term is bounded by $|B_0|\varepsilon^n$;
all other terms follow by the binomial formula. For $m=1$, the
sum is $\rho_0(\varepsilon+\rho_0)^{n-1}$. Hence both the inner
Taylor series and the outer Laurent series converge absolutely
on neighborhoods of the respective physical boundaries.
In particular, the successive translations and Fourier--Gegenbauer
projections are absolutely summable after differentiating once.
Let $E$ be a compact set satisfying
$E\Subset\{|x|<S\}\cap\{|x-\varepsilon e_d|>\rho_0\}$, and put
$r_+:=\sup_E|x|<S$, $r_-:=\inf_E|x-\varepsilon e_d|>\rho_0$.
For each multiindex $\alpha$, the polynomial bounds for the
Gegenbauer harmonics and \eqref{eq:multi-H-rebuild}--
\eqref{eq:multi-B-rebuild} give, for a suitable exponent $M\ge0$,
\begin{equation}\label{eq:multi-rebuild-normal}
 \begin{aligned}
 &\sum_{n\ge m}|A_n|\sup_E|\partial^\alpha R_n^{(m)}|
 +\sum_{j\ge m}|B_j|\sup_E
 |\partial^\alpha I_j^{(m)}(\cdot-\varepsilon e_d)|\\
 &\quad\le C_{\alpha,E}\left[
 \sum_{n\ge1}(n+1)^M\left(\frac{r_+}{S}\right)^n
 +\sum_{j\ge1}(j+1)^M\left(\frac{\rho_0}{r_-}\right)^j
 +1\right]<\infty.
 \end{aligned}
\end{equation}
The finite zero-mode contributions are included in the constant
term in \eqref{eq:multi-rebuild-normal}. To obtain
\eqref{eq:multi-rebuild-normal}, one uses, for each fixed $\alpha$,
solid-harmonic derivative bounds of the form
\begin{equation}\label{eq:multi-harmonic-derivative-bound}
 \begin{aligned}
 \sup_{|x|\le r_+}|\partial^\alpha R_n^{(m)}(x)|
 &\le C_{\alpha,d,r_+}(n+1)^{M_\alpha}
            \max\{r_+,1\}^{|\alpha|}r_+^{\max(n-|\alpha|,0)},\\
 \sup_{|x-\varepsilon e_d|\ge r_-}
 |\partial^\alpha I_n^{(m)}(x-\varepsilon e_d)|
 &\le C_{\alpha,d,r_-}(n+1)^{M_\alpha}
               r_-^{-n-c-|\alpha|}.
 \end{aligned}
\end{equation}
These estimates do not require normalized angular bases.
For $\gamma>0$, factor the generating function as
\[
 (1-2t\cos\theta+t^2)^{-\gamma}
 =(1-te^{i\theta})^{-\gamma}(1-te^{-i\theta})^{-\gamma}.
\]
The triangle inequality for its coefficients gives
$|C_n^\gamma(\cos\theta)|\le C_n^\gamma(1)$, and
$C_n^\gamma(1)=\Gamma(n+2\gamma)/(n!\Gamma(2\gamma))$
has polynomial growth. This controls the undifferentiated
harmonics in \eqref{eq:multi-bases0} and \eqref{eq:associated}.
For a regular harmonic, apply the interior derivative estimate
on a ball of radius comparable to $r_+/(n+1)$ around a point
of $B_{r_+}$, using the bound on
$B_{r_+(1+1/(n+1))}$. The factor
$(1+1/(n+1))^n$ is bounded, and each derivative contributes
only a polynomial factor in $n$ and a power of $r_+^{-1}$.
For a singular harmonic at distance $r\ge r_->0$ from its
center, use a ball of radius $r/(2(n+1))$ and the analogous
bound with $(1-1/(2(n+1)))^{-n-c}$. This proves
\eqref{eq:multi-harmonic-derivative-bound} for every fixed
multiindex; the positive radius $r_+$ may always be enlarged
slightly if necessary.
Only finitely many $n<|\alpha|$ require a separate constant, and
$R_n^{(m)}$ is polynomial. Combining
\eqref{eq:multi-harmonic-derivative-bound} with
\eqref{eq:multi-fixed-decay} and \eqref{eq:multi-B-rebuild}
yields \eqref{eq:multi-rebuild-normal} by comparison with geometric
series. Therefore the reconstructed series and all their derivatives
converge normally across both physical boundaries, and their sum
is a smooth harmonic function on \eqref{eq:domain}.

By the normal convergence just established, all changes of center
in \eqref{eq:reg-shift}--\eqref{eq:multi-HJ} can be evaluated at the
physical boundaries and differentiated term by term. Denote the
angular factors in the $m$ sector by $\psi_j^{(m)}$. For the positive
degrees, the two boundary residuals have the explicit expansions
\begin{align}
 (\partial_\nu u-\lambda u)|_{|x-\varepsilon e_d|=\delta}
 &=\sum_{j\ge\max(1,m)}
 \left[-jH_j\delta^{j-1}+(j+c)B_j\delta^{-j-c-1}
 -\lambda(H_j\delta^j+B_j\delta^{-j-c})\right]
 \psi_j^{(m)},\label{eq:inner-boundary-residual}\\
 (\partial_\nu u-\lambda u)|_{|x|=1}
 &=\sum_{n\ge\max(1,m)}
 \left[nA_n-(n+c)J_n-\lambda(A_n+J_n)\right]
 \psi_n^{(m)}.\label{eq:outer-boundary-residual}
\end{align}
Writing $\mathcal B_{\mathrm i,j}^{(d)}$ and
$\mathcal B_{\mathrm o,n}^{(d)}$ for the brackets in
\eqref{eq:inner-boundary-residual} and
\eqref{eq:outer-boundary-residual}, respectively, one obtains the
factorizations
\begin{equation}\label{eq:multi-residual-factorization}
 \begin{aligned}
 \mathcal B_{\mathrm i,j}^{(d)}
 &=\frac{j+c-\lambda\delta}{\delta^{j+c+1}}
   \left(B_j-\delta^{2j+c}\beta_j^{(d)}H_j\right)=0,\\
 \mathcal B_{\mathrm o,n}^{(d)}
 &=(n-\lambda)\left(A_n-\alpha_n^{(d)}J_n\right)=0.
 \end{aligned}
\end{equation}
The first equality uses \eqref{eq:multi-BC}; the second uses
\eqref{eq:multi-fixed} and \eqref{eq:multi-HJ}.
For $m=0$, the degree-zero residuals are exactly
\eqref{eq:multi-zero} and vanish by \eqref{eq:multi-B0}; for $m=1$
no zero degree occurs. Normal convergence from
\eqref{eq:multi-rebuild-normal} now gives
\begin{equation}\label{eq:multi-boundary-vanishing}
 (\partial_\nu u-\lambda u)|_{\partial\OO}=0.
\end{equation}

For injectivity, suppose the reconstructed $u_m$ vanishes. The
uniqueness of \eqref{eq:multi-split} in \Cref{lem:multi-split}
then gives $A_n=B_n=0$ for every $n\ge1$; in the axial sector,
\eqref{eq:multi-B0} also gives $A_0=B_0=0$.
For the two inverse identities, extraction from a reconstructed
field returns its prescribed $A$ by uniqueness of the Green
splitting. Conversely, if one starts from a Steklov eigenfunction
$u_m$, its extracted $A$ determines every $B_j$ through
\eqref{eq:multi-HJ} and \eqref{eq:multi-BC}, and determines
$A_0,B_0$ through \eqref{eq:multi-B0} when $m=0$.
The reconstructed function therefore equals that original $u_m$.
Finally, for the reconstructed smooth function, integration by parts
gives, for every $\varphi\in H^1(\OO;\mathbb R)$,
\[
 \int_{\OO}\nabla u_m\cdot\nabla\varphi\,dx
 -\lambda\int_{\partial\OO}(\Tr u_m)(\Tr\varphi)\,dS
 =\int_{\partial\OO}(\partial_\nu u_m-\lambda u_m)
                      \Tr\varphi\,dS=0.
\]
Thus it belongs to the weak eigenspace
\eqref{eq:eigenspace-definition}; the reverse passage to a smooth
function uses the regularity stated in \Cref{sec:prelim}.
The forward and backward maps are therefore inverse. Moreover, \eqref{eq:bounded-extraction}
and \eqref{eq:bounded-left-inverse} show that these inverse maps
are bounded for the $H^1(\OO)$ norm on the eigenspace and the
$\ell^2$ norm on the fixed-point kernel. Thus, for $m=0$ or each
fixed $Y_1$-component, \eqref{eq:multi-equiv} is a topological
linear isomorphism and preserves the dimension of the corresponding
real eigenspace.

For the Neumann endpoint $\lambda=0$, the degree-zero equations
\eqref{eq:multi-zero} give $B_0=0$, while $A_0$ is free.
To reconstruct a nonconstant function from a nonzero fixed point
$A\in\ker(\Id-K^{(d,0)}_{0,\varepsilon})$, choose $A_0=0$ and
set
\begin{equation}\label{eq:multi-neumann-reconstruction}
 \begin{gathered}
 A_0=B_0=0,\qquad
 B_j=\delta^{2j+c}\frac{j}{j+c}H_j,\qquad
 A_n=\frac{n+c}{n}J_n,\qquad j,n\ge1,\\
 (\partial_\nu u)|_{\partial\OO}=0,\qquad
 u\equiv\mathrm{constant}\ \Longrightarrow\ A=0.
 \end{gathered}
\end{equation}
The first line follows from the fixed-point equations and
\eqref{eq:multi-HJ}; the boundary identity follows from the
factorizations in \eqref{eq:multi-residual-factorization} with
$\lambda=0$. The final implication follows from the uniqueness of
\eqref{eq:two-center-coefficient-extraction}. Thus every nonzero
fixed point in \eqref{eq:multi-neumann-reconstruction} produces a
\emph{nonconstant} Neumann harmonic function. This construction
uses only \Cref{lem:multi-compact} and the coefficient equations,
not the spectral-radius inequality \eqref{eq:multi-r0}.
For $m=1$, the same reconstruction and positive-degree residual
identities apply at $\lambda=0$, with no degree-zero terms.
A nonzero fixed point would again give a nonconstant Neumann
harmonic function. Green's energy identity therefore implies
\[
 \ker(\Id-K^{(d,m)}_{0,\varepsilon})=\{0\},\qquad m=0,1,
\]
whereas the PDE zero eigenspace consists of constants.
The spectral-radius estimate at $\lambda=0$ follows from
\Cref{lem:multi-neumann}.
\end{proof}

\section{Spectral comparison for \texorpdfstring{$d\ge3$}{d >= 3}}\label{sec:spectral}
We combine the spectral-radius crossing for
\eqref{eq:Kd0}--\eqref{eq:Kd1} with a comparison of angular sectors.

\subsection{Spectral-radius crossing}\label{subsec:radius-cross}
The next lemma verifies the hypotheses of the Perron--Jentzsch theorem for the two coefficient operators.
\begin{lemma}\label{lem:multi-perron}
For $0<\varepsilon<1-\delta$ and $0\le\lambda<1$, both $K^{(d,0)}_{\lambda,\varepsilon}$ and $K^{(d,1)}_{\lambda,\varepsilon}$ are compact, positive and ideal-irreducible. Each has an algebraically simple positive spectral radius and coordinatewise strictly positive left and right eigenvectors. At concentricity,
\begin{equation}\label{eq:multi-diagonal}
 K^{(d,0)}_{\lambda,0}=K^{(d,1)}_{\lambda,0}
 =\operatorname{diag}\left(
 \alpha_n^{(d)}\beta_n^{(d)}\delta^{2n+c}\right)_{n\ge1}.
\end{equation}
\end{lemma}
\begin{proof}
For $n,k\ge1$ and $\varepsilon>0$, the $j=1$ contributions
in \eqref{eq:Kd0}--\eqref{eq:Kd1} give
\begin{equation}\label{eq:perron-entry-lower}
 \begin{aligned}
 (K^{(d,0)}_{\lambda,\varepsilon})_{nk}
 &\ge\alpha_n^{(d)}n\binom{k+c-1}{k-1}
 \delta^{c+2}\beta_1^{(d)}\varepsilon^{n+k-2}>0,\\
 (K^{(d,1)}_{\lambda,\varepsilon})_{nk}
 &\ge\alpha_n^{(d)}\binom{k+c}{k-1}
 \delta^{c+2}\beta_1^{(d)}\varepsilon^{n+k-2}>0.
 \end{aligned}
\end{equation}
The Hilbert--Schmidt property follows from \Cref{lem:multi-compact};
\Cref{lem:prelim-perron} supplies the asserted
left and right eigenvectors and algebraic simplicity.
For $\varepsilon=0$, only terms with $n=k=j$ remain in the two matrices,
so \eqref{eq:Kd0}--\eqref{eq:Kd1} give
\eqref{eq:multi-diagonal}.
\end{proof}

At concentricity, the coefficient operators are diagonal, and the
first diagonal entry determines their spectral radius.
\begin{lemma}\label{lem:multi-concentric-diagonal}
Fix $d\ge3$, $0<\delta<1$, $0\le\lambda<1$, and put $c=d-2$ as in
\eqref{eq:dim-param}. For every $n\ge1$, define
\begin{equation}\label{eq:concentric-highdim-eta}
 \eta_n^{(d)}(\lambda)
 :=\alpha_n^{(d)}(\lambda)\beta_n^{(d)}(\lambda)\delta^{2n+c}.
\end{equation}
Then
\begin{equation}\label{eq:concentric-highdim-order}
 0<\eta_{n+1}^{(d)}(\lambda)<\eta_n^{(d)}(\lambda),
 \qquad
 \spr(K^{(d,0)}_{\lambda,0})
 =\spr(K^{(d,1)}_{\lambda,0})
 =\eta_1^{(d)}(\lambda).
\end{equation}
\end{lemma}
\begin{proof}
The multipliers in \eqref{eq:multi-factors} satisfy
\begin{equation}\label{eq:concentric-product-identity}
 \alpha_n^{(d)}\beta_n^{(d)}
 =1+\lambda(1+\delta)
 \frac{2n+c}{(n-\lambda)(n+c-\lambda\delta)}.
\end{equation}
For real $t\ge1$ and $0\le\lambda<1$, direct differentiation gives
\begin{equation}\label{eq:concentric-rational-derivative}
 \begin{aligned}
 &\frac{d}{dt}\frac{2t+c}{(t-\lambda)(t+c-\lambda\delta)}\\
 &\quad=-\frac{2(t^2-\lambda^2\delta)+2ct+c^2
       +\lambda c(1-\delta)}
       {(t-\lambda)^2(t+c-\lambda\delta)^2}<0.
 \end{aligned}
\end{equation}
Indeed, $t\ge1>\lambda$ and $\delta<1$ imply
$t^2-\lambda^2\delta>0$; the remaining numerator terms in
\eqref{eq:concentric-rational-derivative} are nonnegative.
For $\lambda>0$, \eqref{eq:concentric-product-identity}--
\eqref{eq:concentric-rational-derivative} show that
$n\mapsto\alpha_n^{(d)}\beta_n^{(d)}$ is strictly decreasing;
for $\lambda=0$ the product is identically one.
As $0<\delta^2<1$, \eqref{eq:concentric-highdim-eta} is strictly
decreasing in both cases. The diagonal operator
\eqref{eq:multi-diagonal} has positive entries tending to zero,
so its spectral radius is the first diagonal entry. This is
\eqref{eq:concentric-highdim-order}.
\end{proof}

The Neumann endpoint supplies the lower side of the crossing in all higher dimensions.
\begin{lemma}\label{lem:multi-neumann}
For fixed $d\ge3$, $0<\delta<1$ and all $0\le\varepsilon<1-\delta$,
\begin{equation}\label{eq:multi-r0}
 \spr(K^{(d,0)}_{0,\varepsilon})<1.
\end{equation}
\end{lemma}
\begin{proof}
By \eqref{eq:multi-diagonal} and \eqref{eq:multi-factors},
$\spr(K^{(d,0)}_{0,0})=\delta^d<1$. The family is norm-continuous
by \Cref{lem:multi-compact}. Its spectral radius is positive
because \eqref{eq:diagonal-spr-bound} and \eqref{eq:Kd0} yield
\begin{equation}\label{eq:multi-entry-lower}
 \spr(K^{(d,0)}_{0,\varepsilon})
 \ge (K^{(d,0)}_{0,\varepsilon})_{11}\ge\delta^d.
\end{equation}
The diagonal estimate \eqref{eq:diagonal-spr-bound} supplies the
first inequality in \eqref{eq:multi-entry-lower}. Thus
\eqref{eq:compact-radius-continuity}, applied with
$T_{\varepsilon}=K^{(d,0)}_{0,\varepsilon}$ using
\Cref{lem:multi-compact}, gives continuity in $\varepsilon$.
If \eqref{eq:multi-r0} failed at some $\varepsilon_1>0$, then
$\spr(K^{(d,0)}_{0,0})=\delta^d<1$ and continuity would furnish
$\varepsilon_*\in(0,\varepsilon_1]$ with
$\spr(K^{(d,0)}_{0,\varepsilon_*})=1$.
By \Cref{lem:multi-perron} there is
$0\ne A\in\ker(\Id-K^{(d,0)}_{0,\varepsilon_*})$.
The explicit reconstruction
\eqref{eq:multi-neumann-reconstruction} in
\Cref{prop:multi-equivalence} then gives a nonconstant harmonic
function $u$ with $\partial_\nu u=0$ on
$\partial\Omega_{\delta,\varepsilon_*}^{(d)}$.
Green's identity would give
\begin{equation}\label{eq:multi-neumann-identity}
 \int_{\Omega^{(d)}_{\delta,\varepsilon_*}}|\nabla u|^2\,dx
 =\int_{\partial\Omega^{(d)}_{\delta,\varepsilon_*}}
          u\partial_\nu u\,dS=0,
\end{equation}
so $u$ would be constant on the connected shell, a contradiction. Thus \eqref{eq:multi-r0} holds.
\end{proof}

Write
\begin{equation}\label{eq:R-d}
 \mathcal R_d(\lambda,\varepsilon)
 :=\spr(K^{(d,0)}_{\lambda,\varepsilon}).
\end{equation}
The estimate \eqref{eq:multi-r0} gives the lower endpoint of the
crossing, while the diagonal lower bound
\eqref{eq:multi-blowup} below provides its upper endpoint.
We now identify the unique level-one root of
$\mathcal R_d(\lambda,\varepsilon)=1$.
\begin{proposition}\label{prop:multi-crossing}
For fixed $d\ge3$, $0<\delta<1$ and $0\le\varepsilon<1-\delta$, $\mathcal R_d(\lambda,\varepsilon)$ is continuous and strictly increasing in $\lambda\in[0,1)$, and
\begin{equation}\label{eq:multi-crossing}
 \mathcal R_d(0,\varepsilon)<1,
 \quad\lim_{\lambda\uparrow1}\mathcal R_d(\lambda,\varepsilon)=+\infty,
 \quad\exists!\,\lambda_*(d,\delta,\varepsilon)\in(0,1):
 \mathcal R_d(\lambda_*,\varepsilon)=1.
\end{equation}
For each fixed $0<\lambda<1$, $\varepsilon\mapsto\mathcal R_d(\lambda,\varepsilon)$ is strictly increasing on $[0,1-\delta)$.
\end{proposition}
\begin{proof}
The Robin multipliers and rank-one coefficient satisfy
\begin{equation}\label{eq:multi-lambda-derivative}
 \begin{aligned}
 \frac{d\alpha_n^{(d)}}{d\lambda}
 &=\frac{2n+c}{(n-\lambda)^2}>0,\\
 \frac{d\beta_j^{(d)}}{d\lambda}
 &=\frac{\delta(2j+c)}{(j+c-\lambda\delta)^2}>0,\\
 \frac{d}{d\lambda}\left(\frac{\lambda\delta^{c+1}}{D_\lambda^{(d)}}\right)
 &=\frac{c\delta^{c+1}(1+\delta^{c+1})}{(D_\lambda^{(d)})^2}>0.
 \end{aligned}
\end{equation}
For $\varepsilon>0$ and $0\le\lambda_1<\lambda_2<1$,
\eqref{eq:multi-lambda-derivative} and \eqref{eq:Kd0} imply
$K^{(d,0)}_{\lambda_2,\varepsilon}>K^{(d,0)}_{\lambda_1,\varepsilon}$
entrywise. Let $v_1>0$ be a right Perron vector for the smaller
operator and $w_2>0$ a left Perron vector for the larger one. Then
\begin{equation}\label{eq:multi-lambda-pairing}
 \begin{aligned}
 &[\mathcal R_d(\lambda_2,\varepsilon)-
   \mathcal R_d(\lambda_1,\varepsilon)]\langle w_2,v_1\rangle
 \\&\quad=\langle w_2,
 (K^{(d,0)}_{\lambda_2,\varepsilon}-
 K^{(d,0)}_{\lambda_1,\varepsilon})v_1\rangle>0.
 \end{aligned}
\end{equation}
At $\varepsilon=0$, \eqref{eq:multi-diagonal} and the boundedness
of the Robin multipliers in \eqref{eq:multi-factors} imply
\begin{equation}\label{eq:multi-diagonal-tail}
 0<\alpha_n^{(d)}(\lambda)\beta_n^{(d)}(\lambda)
 \delta^{2n+c}\le C_{d,\delta,\lambda}\delta^{2n+c}
 \xrightarrow[n\to\infty]{}0.
\end{equation}
In fact \Cref{lem:multi-concentric-diagonal} identifies the unique
maximizing index as $n_1=1$, so $e_1$ is the corresponding
nonnegative eigenvector. From \eqref{eq:multi-diagonal},
\begin{equation}\label{eq:multi-diagonal-lambda}
 \mathcal R_d(\lambda_2,0)
 \ge (K^{(d,0)}_{\lambda_2,0})_{n_1n_1}
 >(K^{(d,0)}_{\lambda_1,0})_{n_1n_1}
 =\mathcal R_d(\lambda_1,0).
\end{equation}
Continuity in $\lambda$ follows by applying
\eqref{eq:compact-radius-continuity} to the Hilbert--Schmidt
continuous family in \Cref{lem:multi-compact}. Its positive-radius
hypothesis holds throughout: the $n=k=j=1$ term gives
$\mathcal R_d(\lambda,\varepsilon)
\ge\alpha_1^{(d)}\beta_1^{(d)}\delta^{c+2}>0$. The first inequality in \eqref{eq:multi-crossing} is \Cref{lem:multi-neumann}; the second follows from
\begin{equation}\label{eq:multi-blowup}
 \mathcal R_d(\lambda,\varepsilon)
 \ge (K^{(d,0)}_{\lambda,\varepsilon})_{11}
 \ge \alpha_1^{(d)}\beta_1^{(d)}\delta^{2+c}
 \xrightarrow[\lambda\uparrow1]{}+\infty.
\end{equation}
The intermediate value theorem yields the unique crossing.

To compare displacements, observe in \eqref{eq:Kd0} that $n+k-2j\ge0$. All the sum terms are nondecreasing in $\varepsilon$, while the rank-one term is strictly increasing in every entry whenever $\lambda>0$. Thus, for $\varepsilon_2>\varepsilon_1$,
\begin{equation}\label{eq:multi-estrict}
 K^{(d,0)}_{\lambda,\varepsilon_2}
 >K^{(d,0)}_{\lambda,\varepsilon_1}
 \quad\text{entrywise}.
\end{equation}
More explicitly, the difference of its $(n,k)$ entries is at least
\[
 \frac{\alpha_n^{(d)}\lambda\delta^{c+1}}{D_\lambda^{(d)}}
 \binom{k+c-1}{k}
 (\varepsilon_2^{n+k}-\varepsilon_1^{n+k})>0.
\]
Here $n+k\ge2$, so the strict inequality also holds when
$\varepsilon_1=0$. Let $v_1\ge0$, $v_1\ne0$, be a right
eigenvector of $K^{(d,0)}_{\lambda,\varepsilon_1}$ at its
spectral radius (take $v_1=e_1$ when $\varepsilon_1=0$),
and let $w_2>0$ be a left Perron vector for
$K^{(d,0)}_{\lambda,\varepsilon_2}$. Writing
$r_i=\mathcal R_d(\lambda,\varepsilon_i)$, we obtain
\begin{equation}\label{eq:multi-ecc-pairing}
 (r_2-r_1)\langle w_2,v_1\rangle
 =\langle w_2,
 (K^{(d,0)}_{\lambda,\varepsilon_2}
 -K^{(d,0)}_{\lambda,\varepsilon_1})v_1\rangle>0.
\end{equation}
\end{proof}

\subsection{Weighted comparison of angular sectors}\label{subsec:angular-compare}
The axial matrix \eqref{eq:Kd0} contains a positive rank-one
correction absent from the transverse matrix \eqref{eq:Kd1}.
The weighted coefficient inequality \eqref{eq:weighted-operator}
provides the spectral-radius comparison in \Cref{lem:multi-branch}.
\begin{lemma}\label{lem:multi-branch}
For every integer $d\ge3$, $0<\delta<1$,
$0<\lambda<1$ and $0<\varepsilon<1-\delta$,
\begin{equation}\label{eq:multi-branch-order}
 \spr(K^{(d,1)}_{\lambda,\varepsilon})
 <\spr(K^{(d,0)}_{\lambda,\varepsilon}).
\end{equation}
At $\varepsilon=0$ the two radii coincide.
\end{lemma}
\begin{proof}
Write $K^{(d,0)}=M^{(d,0)}+P^{(d)}$, where
\begin{equation}\label{eq:multi-P}
 (P^{(d)})_{nk}
 =\frac{\alpha_n^{(d)}\lambda\delta^{c+1}}{D_\lambda^{(d)}}
 \bn{k+c-1}{k}\varepsilon^{n+k}>0.
\end{equation}
The two sequences
$(\alpha_n^{(d)}\varepsilon^n)_{n\ge1}$ and
$(\binom{k+c-1}{k}\varepsilon^k)_{k\ge1}$ belong to $\ell^2$:
the first multiplier is bounded, and the binomial coefficient
is $O((k+1)^{c-1})$. Thus $P^{(d)}$ is a bounded rank-one
operator. Since $0\le M^{(d,0)}\le K^{(d,0)}$ entrywise,
$M^{(d,0)}$ satisfies \eqref{eq:multi-HS-bound} as well.
By \eqref{eq:multi-ratio}, for $1\le j\le\min(n,k)$,
\begin{equation}\label{eq:weighted-ratio}
 \frac{n+c}{k+c}
 \frac{\bn{n-1}{j-1}\bn{k+c}{k-j}}
 {\bn nj\bn{k+c-1}{k-j}}
 =\frac{j(n+c)}{n(j+c)}\le1,
\end{equation}
since $j\le n$. The exact deficit is
\[
 1-\frac{j(n+c)}{n(j+c)}=\frac{c(n-j)}{n(j+c)}\ge0.
\]
All binomial factors are positive, and the ratio equals one
precisely when $j=n$. The summands with $j<n$ have a strict
deficit. Summing the nonnegative terms gives
\begin{equation}\label{eq:weighted-operator}
 \frac{n+c}{k+c}(K^{(d,1)}_{\lambda,\varepsilon})_{nk}
 \le (M^{(d,0)}_{\lambda,\varepsilon})_{nk}.
\end{equation}
To specify the domain of the weight, put
\[
 D=\operatorname{diag}(n+c)_{n\ge1},\qquad
 \mathcal D(D)=\left\{y\in\ell^2:
       \sum_{n\ge1}(n+c)^2|y_n|^2<\infty\right\}.
\]
This is a densely defined closed operator, with bounded inverse
$D^{-1}$ on $\ell^2$. The estimate \eqref{eq:multi-HS-bound}
shows that $K^{(d,1)}$ takes every $y\in\ell^2$ into
$\mathcal D(D)$; indeed
\[
 \sum_{n,k\ge1}(n+c)^2|(K^{(d,1)})_{nk}|^2
 \le C^2\left(\sum_{n\ge1}(n+c)^2R^{-2n}\right)
          \left(\sum_{k\ge1}(k+1)^{2c}q^{2k}\right)<\infty.
\]
Hence $DK^{(d,1)}$ is a bounded Hilbert--Schmidt operator.
Let $v=(v_n)_{n\ge1}>0$ be a right Perron eigenvector of
$K^{(d,1)}_{\lambda,\varepsilon}$ with eigenvalue
$r_1=\spr(K^{(d,1)}_{\lambda,\varepsilon})>0$.
Applying \eqref{eq:multi-HS-bound} to
$K^{(d,1)}v=r_1v$ gives
\begin{equation}\label{eq:perron-weighted-decay}
 \begin{aligned}
 |v_n|
 &=r_1^{-1}|(K^{(d,1)}_{\lambda,\varepsilon}v)_n|\\
 &\le\frac{C}{r_1}R^{-n}
 \left(\sum_{k\ge1}(k+1)^{2c}q^{2k}\right)^{1/2}
 \|v\|_{\ell^2},\qquad n\ge1,\\
 \sum_{n\ge1}(n+c)^2|v_n|^2
 &\le C_{d,\delta,\lambda,\varepsilon,v}
 \sum_{n\ge1}(n+c)^2R^{-2n}<\infty.
 \end{aligned}
\end{equation}
Consequently $x_n:=(n+c)v_n$ defines an element
$x\in\ell^2$. Moreover, \eqref{eq:multi-HS-bound} and the
Cauchy--Schwarz inequality yield, for each $n\ge1$,
\begin{equation}\label{eq:weighted-sum-legitimacy}
 \begin{aligned}
 \sum_{k\ge1}|(M^{(d,0)})_{nk}x_k|
 &\le CR^{-n}\sum_{k\ge1}(k+1)^cq^k|x_k|\\
 &\le CR^{-n}\|x\|_{\ell^2}
 \left(\sum_{k\ge1}(k+1)^{2c}q^{2k}\right)^{1/2}<\infty.
 \end{aligned}
\end{equation}
For each fixed $n$, the other side of the comparison is also
absolutely summable, since
\[
 \sum_{k\ge1}\frac{n+c}{k+c}(K^{(d,1)})_{nk}x_k
 =(n+c)(K^{(d,1)}v)_n=r_1x_n.
\]
One may first sum \eqref{eq:weighted-operator} over $1\le k\le N$
and then let $N\to\infty$ in that coordinate. Thus multiplying
by $(k+c)v_k=x_k$ and summing gives
\begin{equation}\label{eq:multi-subeigen}
 K^{(d,0)}x\ge r_1x+P^{(d)}x>r_1x
 \quad\text{coordinatewise}.
\end{equation}
Let $w>0$ be a left Perron eigenvector of $K^{(d,0)}$ with
eigenvalue $r_0$. Every term in the following double sum is
nonnegative, and the pairing estimate in the proof of
\Cref{lem:prelim-radius-order} gives
\[
 \sum_{n,k\ge1}w_n\bigl((M^{(d,0)})_{nk}+(P^{(d)})_{nk}\bigr)x_k
 =\langle w,K^{(d,0)}x\rangle
 \le\|w\|_{\ell^2}\|K^{(d,0)}\|\|x\|_{\ell^2}<\infty.
\]
Also $0<\langle w,x\rangle<\infty$. We may therefore pair
\eqref{eq:multi-subeigen} with $w$ and factor the rank-one
series without any conditional rearrangement.
Since $0<\varepsilon<1-\delta$, the
rank-one formula \eqref{eq:multi-P} gives the explicit
strictly positive pairing
\begin{equation}\label{eq:rank-one-strict-pairing}
 \begin{aligned}
 \langle w,P^{(d)}x\rangle_{\ell^2}
 &=\frac{\lambda\delta^{c+1}}{D_\lambda^{(d)}}
 \left(\sum_{n\ge1}w_n\alpha_n^{(d)}\varepsilon^n\right)
 \left(\sum_{k\ge1}x_k\binom{k+c-1}{k}\varepsilon^k\right)>0.
 \end{aligned}
\end{equation}
Both scalar series in \eqref{eq:rank-one-strict-pairing}
converge absolutely by Cauchy--Schwarz, $w,x\in\ell^2$,
the boundedness of $\alpha_n^{(d)}$ for fixed $\lambda<1$,
and $\varepsilon<1$. Thus
\begin{equation}\label{eq:multi-pairing}
 (r_0-r_1)\langle w,x\rangle
 =\langle w,K^{(d,0)}x-r_1x\rangle
 \ge\langle w,P^{(d)}x\rangle>0.
\end{equation}
Hence \eqref{eq:multi-pairing} proves
\eqref{eq:multi-branch-order}; the concentric identity is
\eqref{eq:multi-diagonal}.
\end{proof}

\begin{remark}\label{rem:unbounded-weight}
The diagonal multiplier $D=\operatorname{diag}(n+d-2)_{n\ge1}$
is unbounded on $\ell^2$. Accordingly, the inequality
\eqref{eq:weighted-operator} is used on the admissible vector
$x=Dv\in\ell^2$, where $v$ is the transverse Perron vector.
The smoothing estimate \eqref{eq:perron-weighted-decay} establishes
$v\in\mathcal D(D)$, and \eqref{eq:weighted-sum-legitimacy}
justifies the coordinatewise sums. The comparison and the pairing
in \eqref{eq:multi-subeigen}--\eqref{eq:multi-pairing} therefore do
not rely on a bounded similarity transformation.
\end{remark}

It remains to eliminate all transverse angular sectors of degree at least two. Their comparison is independent of the infinite-matrix method.
\begin{lemma}\label{lem:higher-angular}
Fix $d\ge3$, $0<\delta<1$, and
$0\le\varepsilon<1-\delta$. For each integer $m\ge1$, let
$\tau_m^{(d)}(\delta,\varepsilon)$ be the smallest Steklov
eigenvalue of \eqref{eq:steklov} restricted to a fixed transverse
spherical-harmonic sector of degree $m$ on $\mathbb S^{d-2}$.
Writing $\tau_m=\tau_m^{(d)}(\delta,\varepsilon)$ within
this lemma, one has
\begin{equation}\label{eq:higher-angular}
 \tau_1<\tau_m\qquad(m\ge2).
\end{equation}
\end{lemma}
\begin{proof}
Write $x=(s\omega,\xi)$ with $s>0$, $\omega\in\mathbb S^{d-2}$, and let $\mathcal D$ be the meridional cross-section of \eqref{eq:domain}. Let $\partial\mathcal D_{\rm sph}$ denote the two circular arcs
of $\partial\mathcal D$ corresponding to the physical inner and outer
spheres; the rotation axis $s=0$ is not part of this physical boundary.
Choose $Y_m$ with $\|Y_m\|_{L^2(\mathbb S^{d-2})}=1$ and
$-\Delta_{\mathbb S^{d-2}}Y_m=m(m+d-3)Y_m$.
For $u=v(s,\xi)Y_m(\omega)$, its Rayleigh quotient is
\begin{equation}\label{eq:angular-quotient}
 \mathcal Q_m^{\mathrm{ang}}[v]=
 \frac{\displaystyle\int_{\mathcal D}
 \left(|v_s|^2+|v_\xi|^2+\frac{m(m+d-3)}{s^2}|v|^2\right)
 s^{d-2}\,ds\,d\xi}
 {\displaystyle\int_{\partial\mathcal D_{\rm sph}}
 |v|^2s^{d-2}\,d\ell}.
\end{equation}
For $m\ge1$, $s\le1$ and the angular eigenvalue satisfies
$\kappa_m:=m(m+d-3)\ge d-2>0$. Consequently,
\begin{equation}\label{eq:angular-coercivity}
 \int_{\mathcal D}|v|^2s^{d-2}\,ds\,d\xi
 \le\int_{\mathcal D}|v|^2s^{d-4}\,ds\,d\xi
 \le\frac1{\kappa_m}\int_{\mathcal D}
 \frac{\kappa_m}{s^2}|v|^2s^{d-2}\,ds\,d\xi.
\end{equation}
Thus the numerator of \eqref{eq:angular-quotient} controls the
$H^1$ norm in the closed $m$-sector. The compact trace map
\eqref{eq:trace-compact} and the direct method give a minimizer
$v_m$ satisfying
$\int_{\partial\mathcal D_{\rm sph}}|v_m|^2s^{d-2}\,d\ell=1$.
The denominator constraint is preserved under a weakly convergent
minimizing sequence because the trace converges strongly in $L^2$;
the numerator is weakly lower semicontinuous. The Euler--Lagrange
equation initially holds for tests in the fixed harmonic sector.
For an arbitrary $H^1$ test, apply its sector projection $P$ and
use the form identities following
\eqref{eq:sector-projection-commutators}; the equation then holds
for the original test as well. Thus the minimizer is an eigenfunction
of the full weak Steklov problem, so this variational minimum is
indeed $\tau_m$. Since $m\ge1$, its boundary mean vanishes by
angular orthogonality, and its positive energy gives $\tau_m>0$.
In particular,
\begin{equation}\label{eq:angular-finite-energy}
 \int_\mathcal D(|\partial_s v_m|^2+|\partial_\xi v_m|^2)
       s^{d-2}\,ds\,d\xi<\infty,\qquad
 \int_\mathcal D|v_m|^2s^{d-4}\,ds\,d\xi<\infty.
\end{equation}
The second integral is the angular part of its finite $m$-energy.
Replacing $Y_m$ by a normalized first-degree transverse harmonic
$Y_1$ leaves the radial-gradient integral and boundary trace norm
unchanged while decreasing the angular coefficient from
$m(m+d-3)$ to $d-2$. Hence $v_mY_1\in H^1(\OO)$, even at the rotation
axis, by the following cutoff argument. Choose a smooth
$\chi_\eta(s)$ that is zero for $s\le\eta$, one for
$s\ge2\eta$, and satisfies $|\chi_\eta'|\le C/\eta$.
The function $\chi_\eta v_mY_1$ belongs to $H^1(\OO)$,
and the only additional gradient term is bounded by
\[
 \int_{\eta<s<2\eta}|\chi_\eta'v_m|^2s^{d-2}\,ds\,d\xi
 \le C\int_{\eta<s<2\eta}|v_m|^2s^{d-4}\,ds\,d\xi
 \longrightarrow0.
\]
All other discarded energy integrals tend to zero by
\eqref{eq:angular-finite-energy}. These functions therefore
converge in $H^1$ to $v_mY_1$, and trace continuity identifies
its physical-boundary norm with the denominator in
\eqref{eq:angular-quotient}. This also covers $d=3$, where
the second weight is $s^{-1}$. Since
$m(m+d-3)>d-2$ for $m\ge2$,
\begin{equation}\label{eq:angular-diff}
 \begin{aligned}
 \tau_1\le \mathcal Q_1^{\mathrm{ang}}[v_m]
 &=\mathcal Q_m^{\mathrm{ang}}[v_m]
 -\bigl(m(m+d-3)-(d-2)\bigr)
 \frac{\displaystyle\int_{\mathcal D}|v_m|^2s^{d-4}\,ds\,d\xi}
 {\displaystyle\int_{\partial\mathcal D_{\rm sph}}|v_m|^2s^{d-2}\,d\ell}\\
 &<\mathcal Q_m^{\mathrm{ang}}[v_m]=\tau_m.
 \end{aligned}
\end{equation}
\end{proof}

\subsection{Identification of the first eigenvalue}\label{subsec:proof-main}
We can now identify the crossing in \eqref{eq:multi-crossing} with the \emph{first} positive Steklov eigenvalue of the full $d$-dimensional shell, not just an eigenvalue of an axial subproblem.
\begin{proposition}\label{prop:multi-identify}
For every $d\ge3$ and admissible $(\delta,\varepsilon)$,
\begin{equation}\label{eq:multi-identify}
 \lambda_*(d,\delta,\varepsilon)
 =\sigma_1^{(d)}(\delta,\varepsilon).
\end{equation}
For $0<\varepsilon<1-\delta$, the eigenspace
\eqref{eq:eigenspace-definition} at the eigenvalue
\eqref{eq:multi-identify} is one-dimensional and lies in the
axial sector $m=0$.
\end{proposition}
\begin{proof}
At the level-one crossing, \Cref{lem:multi-perron} and the diagonal endpoint \eqref{eq:multi-diagonal} yield a nonzero fixed point of $K^{(d,0)}$. By \Cref{prop:multi-equivalence} it reconstructs a Steklov eigenfunction, so $\lambda_*$ is a positive eigenvalue. For $0<\lambda<\lambda_*$, \Cref{prop:multi-crossing,lem:multi-branch} give
\begin{equation}\label{eq:multi-exclude}
 \spr(K^{(d,1)}_{\lambda,\varepsilon})
 \le \spr(K^{(d,0)}_{\lambda,\varepsilon})<1.
\end{equation}
At $\lambda=\lambda_*$, the transverse radius is strictly
below one if $\varepsilon>0$, by \eqref{eq:multi-branch-order}.
If $\varepsilon=0$, the two operators agree for every
$0\le\lambda<1$, by \eqref{eq:multi-diagonal}; their common
radius equals one at $\lambda_*$ and is below one for
$0<\lambda<\lambda_*$. By \Cref{prop:multi-equivalence}, this excludes
$m=1$ eigenvalues on $(0,\lambda_*]$ when $\varepsilon>0$ and
on $(0,\lambda_*)$ when $\varepsilon=0$. Define $\tau_m$ as in
\Cref{lem:higher-angular}. Then
\begin{equation}\label{eq:multi-sector-exclusion}
 \begin{cases}
 \tau_1>\lambda_*,&\varepsilon>0,\\
 \tau_1\ge\lambda_*,&\varepsilon=0,
 \end{cases}
 \qquad
 \tau_m>\tau_1\quad(m\ge2).
\end{equation}
To make the exclusion of every smaller full eigenvalue explicit,
suppose $0<\lambda<\lambda_*$ were such an eigenvalue and take
a nonzero real eigenfunction $u$. The bounded angular projections
following \eqref{eq:sector-projection-commutators} provide a nonzero
fixed harmonic component with the same eigenvalue. An axial
component contradicts \Cref{prop:multi-equivalence} and
\eqref{eq:multi-exclude}; a degree-one transverse component does
the same. A component of degree $m\ge2$ would give
$\lambda\ge\tau_m>\tau_1\ge\lambda_*$ by the attained sectorial
variational minimum, again a contradiction.
For $\lambda>0$, testing \eqref{eq:eigenspace-definition}
with $\varphi\equiv1$ gives
\[
 0=\lambda\int_{\partial\OO}\Tr u\,dS,
 \qquad u\equiv C\ \Longrightarrow\ C=0.
\]
Thus the axial crossing is the smallest positive eigenvalue of
the full problem. At $\varepsilon=0$, the diagonal fixed point
$e_1$ also reconstructs in every fixed degree-one transverse
component, so the non-strict inequality in
\eqref{eq:multi-sector-exclusion} sharpens to $\tau_1=\lambda_*$.
There are $d-1$ independent transverse directions; the scalar
coefficient kernel represents each one separately.
For $\varepsilon>0$, the Perron eigenspace of
$K^{(d,0)}_{\lambda_*,\varepsilon}$ is one-dimensional by
\cite[Theorem~1.5]{ChangWangWu2020}. The coefficient--PDE
isomorphism in \Cref{prop:multi-equivalence} preserves this
dimension. By \eqref{eq:multi-sector-exclusion}, neither a
fixed $Y_1$-component nor any degree $m\ge2$ contributes to
this eigenspace. Thus
\begin{equation}\label{eq:multi-eigenspace-isomorphism}
 \dim_{\mathbb R}\mathscr E_{\lambda_*}^{(d)}(\delta,\varepsilon)
 =\dim_{\mathbb R}\ker(\Id-K^{(d,0)}_{\lambda_*,\varepsilon})=1,
 \qquad \varepsilon>0.
\end{equation}
\end{proof}

We complete the proof for $d\ge3$ and compute the concentric
multiplicity by separation of variables.
\begin{proof}[Proof of \Cref{thm:main} for $d\ge3$]
Let $\lambda_i=\sigma_1^{(d)}(\delta,\varepsilon_i)$, $i=1,2$. By \eqref{eq:multi-identify} and \eqref{eq:multi-crossing},
\begin{equation}\label{eq:multi-roots}
 \mathcal R_d(\lambda_1,\varepsilon_1)=1,
 \qquad \mathcal R_d(\lambda_2,\varepsilon_2)=1.
\end{equation}
By \Cref{prop:multi-crossing}, the comparison at the first root
$\lambda_1$ in \eqref{eq:multi-roots} gives
\begin{equation}\label{eq:multi-monotone-chain}
 1=\mathcal R_d(\lambda_1,\varepsilon_1)
 <\mathcal R_d(\lambda_1,\varepsilon_2),\qquad
 1=\mathcal R_d(\lambda_2,\varepsilon_2).
\end{equation}
Since $\lambda\mapsto\mathcal R_d(\lambda,\varepsilon_2)$
is strictly increasing, \eqref{eq:multi-monotone-chain} gives
$\lambda_2<\lambda_1$, as required in
\eqref{eq:main-monotonicity}.

At $\varepsilon=0$, for $0<\lambda<1$, separation into spherical
harmonics reduces all positive angular degrees to the diagonal
coefficient equations
\begin{equation}\label{eq:concentric-highdim-mode-equation}
 \eta_n^{(d)}(\lambda)=1,\qquad n\ge1,
\end{equation}
with $\eta_n^{(d)}$ defined in
\eqref{eq:concentric-highdim-eta}.
By \Cref{lem:multi-concentric-diagonal}, for the first crossing
$\lambda_*:=\lambda_*(d,\delta,0)$ from
\eqref{eq:multi-crossing},
\begin{equation}\label{eq:concentric-highdim-mode-exclusion}
 \begin{aligned}
 \eta_1^{(d)}(\lambda_*)&=1,\\
 0<\lambda\le\lambda_*,\quad n\ge2
 &\quad\Longrightarrow\quad
 \eta_n^{(d)}(\lambda)
 <\eta_1^{(d)}(\lambda)\le1.
 \end{aligned}
\end{equation}
Thus the only nonconstant angular degree at the first crossing
is $n=1$.
The degree-zero harmonic functions are $u(r)=A+Br^{-c}$.
For a nonconstant radial eigenfunction $B\ne0$ and $\lambda>0$, the two radial boundary equations from
\eqref{eq:steklov} give
\begin{equation}\label{eq:radial-zero-system}
 \begin{aligned}
 -cB&=\lambda(A+B),\\
 cB\delta^{-c-1}&=\lambda(A+B\delta^{-c}),
 \end{aligned}
 \qquad\Longrightarrow\qquad
 \lambda(\delta^{-c}-1)=c(\delta^{-c-1}+1).
\end{equation}
Besides the constant zero mode, the only radial Steklov eigenvalue is
\begin{equation}\label{eq:radial-zero-d}
 \Lambda_{0,+}
 =\frac{c(1+\delta^{c+1})}{\delta(1-\delta^c)}>1.
\end{equation}
By \eqref{eq:range} and \eqref{eq:radial-zero-d}, the degree-zero positive radial eigenvalue lies strictly above the interval $(0,1)$ containing the first eigenvalue. Together with \eqref{eq:concentric-highdim-mode-exclusion}, this shows that the first positive eigenspace has angular degree $n=1$. Its angular dimension is $d$, and the two radial boundary equations determine the coefficient ratio uniquely. Hence its multiplicity is $d$. Simplicity for $\varepsilon>0$ follows from \Cref{prop:multi-identify}.
\end{proof}

\section{The planar case}\label{sec:plane}

For $d=2$, the logarithmic fundamental solution leads to a
holomorphic coefficient representation. Identify $(x_1,x_2)$ with $z=x_2+i x_1$, so that the displacement $\varepsilon e_2$ corresponds to the real point $\varepsilon$. Reflection $z\mapsto\bar z$ fixes the domain and splits the spectrum into even and odd real sectors.

\subsection{Holomorphic expansions and coefficient equations}\label{subsec:plane-reduce}
The harmonic decomposition below is the planar analogue of
\Cref{lem:multi-split}; the logarithmic term represents the
degree-zero multipole. See the harmonic decomposition and annular-region theory in
\cite[Chs.~10--11, pp.~191--221]{AxlerBourdonRamey}. Normalization of both holomorphic parts makes it unique.
\begin{lemma}\label{lem:planar-split}
Every real harmonic $u\in C^2(\Omega^{(2)}_{\delta,\varepsilon})$ has a unique representation
\begin{equation}\label{eq:planar-split}
 u(z)=C+b\log|z-\varepsilon|+\operatorname{Re}\left(
 \sum_{n\ge1}A_nz^n+\sum_{j\ge1}\frac{B_j}{(z-\varepsilon)^j}\right),\quad C,b\in\RR,
\end{equation}
where the regular and singular series converge normally in $|z|<1$ and $|z-\varepsilon|>\delta$, respectively. If $u$ is a smooth Steklov eigenfunction, $A=(A_n)_{n\ge1}\in\ell^2(\CC)$.
\end{lemma}
\begin{proof}
The Wirtinger derivative $g=2\partial_z u=u_{x_2}-iu_{x_1}$ is holomorphic. Cauchy's formula on two auxiliary circles $|z|=r_{\mathrm o}<1$ and $|z-\varepsilon|=r_{\mathrm i}>\delta$ with $\varepsilon+r_{\mathrm i}<r_{\mathrm o}$ yields
\begin{equation}\label{eq:planar-h}
 g(z)=\sum_{n\ge0}p_nz^n+
       \sum_{n\ge0}\frac{q_n}{(z-\varepsilon)^{n+1}},\qquad
 \begin{cases}
 p_n=\dfrac1{2\pi i}\displaystyle\oint_{|\zeta|=r_{\mathrm o}}\dfrac{g(\zeta)}{\zeta^{n+1}}\,d\zeta,\\[4pt]
 q_n=\dfrac1{2\pi i}\displaystyle\oint_{|\zeta-\varepsilon|=r_{\mathrm i}}
 g(\zeta)(\zeta-\varepsilon)^n\,d\zeta.
 \end{cases}
\end{equation}
Both contours are counterclockwise; the negative sign in the inner boundary of the doubly connected Cauchy formula cancels the negative sign in the inner kernel expansion. Set $M_{\mathrm o}:=\sup_{|z|=r_{\mathrm o}}|g(z)|$ and
$M_{\mathrm i}:=\sup_{|z-\varepsilon|=r_{\mathrm i}}|g(z)|$.
Cauchy estimates give
\begin{equation}\label{eq:planar-h-est}
 |p_n|\le M_{\mathrm o}r_{\mathrm o}^{-n},\qquad
 |q_n|\le M_{\mathrm i}r_{\mathrm i}^{n+1}.
\end{equation}
Deforming the auxiliary circles shows that the coefficients are independent of their chosen radii. For every closed piecewise-$C^1$ path $\gamma$ in the shell,
\begin{equation}\label{eq:period-real}
 \operatorname{Re}\oint_\gamma g(z)\,dz
 =\oint_\gamma du=0.
\end{equation}
In particular, integration around the inner circle gives
$\operatorname{Re}(2\pi i q_0)=0$, hence $q_0\in\mathbb R$.
The single-valued real harmonic function therefore contributes
a real logarithmic coefficient, with no multivalued argument
term. The remaining Cauchy modes have the single-valued
primitives used below. Define
\begin{equation}\label{eq:planar-coeff-deriv}
 b=q_0,\qquad A_n=p_{n-1}/n,\qquad B_j=-q_j/j.
\end{equation}
Let $F(z)=\sum_{n\ge1}A_nz^n$ and
$G(z)=\sum_{j\ge1}B_j(z-\varepsilon)^{-j}$,
with $A_n,B_j$ from \eqref{eq:planar-coeff-deriv}. Since the two
series converge locally normally, termwise differentiation gives
\begin{equation}\label{eq:planar-primitive-check}
 2\partial_z\left[b\log|z-\varepsilon|+
       \operatorname{Re}(F(z)+G(z))\right]
 =\frac b{z-\varepsilon}+F'(z)+G'(z)=g(z).
\end{equation}
The difference between $u$ and the bracketed function in
\eqref{eq:planar-primitive-check} has zero gradient on the connected
shell; it is a real constant $C$. The uniqueness of the Cauchy
coefficients in \eqref{eq:planar-h}, together with the value of $C$,
proves uniqueness in \eqref{eq:planar-split}. In this uniqueness
statement the holomorphic normalizations are $F(0)=0$ and
$G(z)\to0$ at infinity; neither holomorphic part carries an
arbitrary imaginary constant. Reflection of the unique
representation shows that the even sector has $A_n,B_j\in\RR$,
whereas the odd sector has $A_n,B_j\in i\RR$ and $C=b=0$.
If $u$ is smooth on the closed shell, its singular derivative is holomorphic across the outer circle. Thus the regular derivative $g_+(z):=\sum_{n\ge0}p_nz^n$ has an $L^2$ outer trace, and Parseval gives
\begin{equation}\label{eq:planar-l2}
 \sum_{n\ge1}|A_n|^2\le\sum_{n\ge0}|p_n|^2
 =\frac1{2\pi}\int_0^{2\pi}|g_+(e^{i\theta})|^2\,d\theta<\infty.
\end{equation}
\end{proof}

The two boundary circles are expanded in their respective center coordinates. We retain the conjugate coefficients explicitly, because the harmonic function is the real part of a complex series.
\begin{proposition}\label{prop:planar-matrix}
Fix $0<\delta<1$ and $0\le\varepsilon<1-\delta$.
For $0\le\lambda<1$, put
\begin{equation}\label{eq:planar-factors}
 \alpha_n^{(2)}=\frac{n+\lambda}{n-\lambda},\qquad
 \beta_j^{(2)}=\frac{j+\lambda\delta}{j-\lambda\delta},\qquad
 D_\lambda^{(2)}=1+\delta^{-1}-\lambda\log(\delta^{-1})>0.
\end{equation}
Indeed, for $t=\delta^{-1}>1$ and $0\le\lambda<1$,
\begin{equation}\label{eq:planar-D-lower}
 D_\lambda^{(2)}=1+t-\lambda\log t
 >1+t-(t-1)=2.
\end{equation}
For $0<\lambda<1$, the boundary equations of
\eqref{eq:steklov} imply
\begin{equation}\label{eq:planar-fixed}
 A^{\mathrm R}=K^{(2,0)}_{\lambda,\varepsilon}A^{\mathrm R},\qquad
 A^{\mathrm I}=K^{(2,1)}_{\lambda,\varepsilon}A^{\mathrm I},
\end{equation}
where $A^{\mathrm R}=\operatorname{Re}A$, $A^{\mathrm I}=\operatorname{Im}A$ and, for $n,k\ge1$,
\begin{align}
 (K^{(2,1)}_{\lambda,\varepsilon})_{nk}
 &=\alpha_n^{(2)}\sum_{j=1}^{\min(n,k)}
   \bn{n-1}{j-1}\bn{k}{j}\delta^{2j}\beta_j^{(2)}
   \varepsilon^{n+k-2j},\label{eq:K21}\\
 (K^{(2,0)}_{\lambda,\varepsilon})_{nk}
 &=(K^{(2,1)}_{\lambda,\varepsilon})_{nk}
   +\frac{\alpha_n^{(2)}\lambda}{nD_\lambda^{(2)}}
      \varepsilon^{n+k}.
 \label{eq:K20}
\end{align}
The additive term in \eqref{eq:K20} has rank one if
$\lambda\varepsilon>0$ and equals the zero operator otherwise.
The reverse implication from \eqref{eq:planar-fixed} to
\eqref{eq:steklov} is established in
\Cref{prop:planar-equivalence}.
\end{proposition}
\begin{proof}
Let $F(z)=\sum_{n\ge1}A_nz^n$ and $G(z)=\sum_{j\ge1}B_j(z-\varepsilon)^{-j}$, as in \eqref{eq:planar-split}. For $|w|<1-\varepsilon$ and $|z|>\varepsilon$, Taylor's theorem
and the binomial series give
\begin{equation}\label{eq:planar-shift-expansions}
 \begin{aligned}
 F(\varepsilon+w)&=F(\varepsilon)+\sum_{j\ge1}H_jw^j,\\
 (z-\varepsilon)^{-j}
 &=\sum_{n\ge j}\binom{n-1}{j-1}\varepsilon^{n-j}z^{-n}.
 \end{aligned}
\end{equation}
Consequently the translation coefficients are
\begin{equation}\label{eq:planar-shifts}
 H_j=\sum_{k\ge j}\bn{k}{j}\varepsilon^{k-j}A_k,
 \qquad
 J_n=\sum_{j=1}^n\bn{n-1}{j-1}\varepsilon^{n-j}B_j.
\end{equation}
For a $2\pi$-periodic function $f$, define its $j$th complex
Fourier coefficient by
\begin{equation}\label{eq:planar-fourier-notation}
 [e^{ij\phi}]f=\frac1{2\pi}\int_0^{2\pi}f(\phi)e^{-ij\phi}\,d\phi.
\end{equation}
For $j\ge1$ and $a\in\CC$, this convention gives
\[
 [e^{ij\phi}]\bigl(2\operatorname{Re}(ae^{ij\phi})\bigr)=a,
 \qquad
 [e^{ij\phi}]\bigl(2\operatorname{Re}(ae^{-ij\phi})\bigr)
 =\overline a.
\]
By contrast, the zero coefficient of a real constant $2C$ is
$2C$ itself. These identities fix both the conjugation of
singular coefficients and the factor two in the zero residuals.
Write $z=\varepsilon+\rho e^{i\phi}$ on the inner circle and
$z=re^{i\theta}$ on the outer circle. The positive Fourier-mode
coefficients of $2u$ and $2\partial_\rho u$ at radius $\rho$ are
\begin{equation}\label{eq:planar-inner-fourier}
 \begin{aligned}
 [e^{ij\phi}](2u)&=H_j\rho^j+\overline{B_j}\rho^{-j},\\
 [e^{ij\phi}](2\partial_\rho u)
 &=jH_j\rho^{j-1}-j\overline{B_j}\rho^{-j-1}.
 \end{aligned}
\end{equation}
On the outer circle, use the normally convergent expansion
\[
 \log|z-\varepsilon|=
 \log r-\operatorname{Re}\sum_{n\ge1}\frac{\varepsilon^n}{nz^n}
 \qquad(r>\varepsilon).
\]
The positive coefficients of $2u$ and $2\partial_r u$ are therefore
\begin{equation}\label{eq:planar-outer-fourier}
 \begin{aligned}
 [e^{in\theta}](2u)&=A_nr^n+
 \left(\overline{J_n}-\frac{b\varepsilon^n}{n}\right)r^{-n},\\
 [e^{in\theta}](2\partial_r u)&=nA_nr^{n-1}-
 n\left(\overline{J_n}-\frac{b\varepsilon^n}{n}\right)r^{-n-1}.
 \end{aligned}
\end{equation}
Apply the inner normal $-\partial_\rho$ at $\rho=\delta$ and
outer normal $\partial_r$ at $r=1$. Equations
\eqref{eq:planar-inner-fourier}--\eqref{eq:planar-outer-fourier}
give
\begin{equation}\label{eq:planar-BC}
 \overline{B_j}=\delta^{2j}\beta_j^{(2)}H_j,
 \qquad
 A_n=\alpha_n^{(2)}\left(\overline{J_n}-\frac{b\varepsilon^n}{n}\right).
\end{equation}
The identities \eqref{eq:planar-inner-fourier} and
\eqref{eq:planar-outer-fourier} hold at interior radii. Smooth
traces from \Cref{lem:planar-split} justify the modewise limits
\begin{equation}\label{eq:planar-Fourier-boundary-limits}
 \begin{aligned}
 [e^{ij\phi}]\bigl(2\partial_\rho u(\delta,\cdot)\bigr)
 &=\lim_{\rho\downarrow\delta}[e^{ij\phi}]
   \bigl(2\partial_\rho u(\rho,\cdot)\bigr),\\
 [e^{in\theta}]\bigl(2\partial_r u(1,\cdot)\bigr)
 &=\lim_{r\uparrow1}[e^{in\theta}]
   \bigl(2\partial_r u(r,\cdot)\bigr).
 \end{aligned}
\end{equation}
Together with the corresponding limits for $u$, this gives
\eqref{eq:planar-BC}. Their degree-zero projections give
\begin{equation}\label{eq:planar-zeromode}
 b=\lambda C,\qquad
 -b/\delta=\lambda\bigl(C+b\log\delta+\operatorname{Re}F(\varepsilon)\bigr),
 \quad
 C=-\frac{\operatorname{Re}F(\varepsilon)}{D_\lambda^{(2)}}
 \quad(\lambda>0).
\end{equation}
For $0<\lambda<1$, the zero-mode equations
\eqref{eq:planar-zeromode} yield
\[
 b=-\frac{\lambda}{D_\lambda^{(2)}}
 \operatorname{Re}F(\varepsilon),\qquad
 \operatorname{Re}F(\varepsilon)
 =\sum_{k\ge1}\varepsilon^k\operatorname{Re}A_k.
\]
Moreover, the translated singular coefficients satisfy
\begin{equation}\label{eq:planar-J-bar}
 \overline{J_n}
 =\sum_{k\ge1}\left[
 \sum_{j=1}^{\min(n,k)}
 \binom{n-1}{j-1}\binom{k}{j}
 \delta^{2j}\beta_j^{(2)}\varepsilon^{n+k-2j}
 \right]A_k.
\end{equation}
Insert \eqref{eq:planar-J-bar} into the second relation of
\eqref{eq:planar-BC}; then
\begin{equation}\label{eq:planar-final-coeff}
 A_n=\sum_{k\ge1}
 (K^{(2,1)}_{\lambda,\varepsilon})_{nk}A_k
 +\frac{\alpha_n^{(2)}\lambda}{nD_\lambda^{(2)}}
 \sum_{k\ge1}\varepsilon^{n+k}\operatorname{Re}A_k.
\end{equation}
The matrices in \eqref{eq:planar-final-coeff} are real. Taking
real and imaginary parts gives \eqref{eq:planar-fixed} and
\eqref{eq:K21}--\eqref{eq:K20}. The convergence and reverse
implication are proved in
\Cref{lem:planar-compact,prop:planar-equivalence}.
At $\lambda=0$, the rank-one term in \eqref{eq:K20} vanishes, so
$K^{(2,0)}_{0,\varepsilon}=K^{(2,1)}_{0,\varepsilon}$; the
constant mode in \eqref{eq:planar-zeromode} then has $b=0$ and
arbitrary $C$.
\end{proof}

\subsection{Compactness and reconstruction}\label{subsec:plane-reconstruct}
The following bounds establish Hilbert--Schmidt compactness and the normal convergence required to reconstruct the boundary conditions.
\begin{lemma}\label{lem:planar-compact}
For fixed $0<\delta<1$, $0\le\lambda<1$,
and $0\le\varepsilon<1-\delta$, the matrices in \eqref{eq:K21}--\eqref{eq:K20} define compact positive operators on $\ell^2(\NN;\RR)$, continuously depending on $(\lambda,\varepsilon)$ in Hilbert--Schmidt norm. Locally uniformly in these parameters there exist $R>1$, $q\in(0,1)$ and $C<\infty$ such that
\begin{equation}\label{eq:planar-hs}
 0\le (K^{(2,m)}_{\lambda,\varepsilon})_{nk}
 \le CR^{-n}q^k,\qquad m=0,1.
\end{equation}
Every fixed point has $|A_n|\le C_AR^{-n}$.
\end{lemma}
\begin{proof}
Recall the shared majorant $T_\varepsilon$ from \eqref{eq:Tmap}. Since $\delta<1-\varepsilon$, $T_\varepsilon(1)<1$. Choose $R>1$ sufficiently close to one that $\varepsilon R<1$ and $q:=T_\varepsilon(R)<1$. The nonnegative majorant
\begin{equation}\label{eq:planar-majorant}
 U^{(2)}_{nk}=\sum_{j=1}^{\min(n,k)}\bn nj\bn kj
 \delta^{2j}\varepsilon^{n+k-2j}
\end{equation}
satisfies, by the binomial series \eqref{eq:binom-series} and the
binomial theorem, the generating identity
\begin{equation}\label{eq:planar-gf}
 \sum_{n\ge1}U^{(2)}_{nk}z^n
 =\frac{T_\varepsilon(z)^k-\varepsilon^k}{1-\varepsilon z},
 \qquad z\ge0,\quad \varepsilon z<1.
\end{equation}
For $z=R$, positivity of the power series in
\eqref{eq:planar-gf} gives
\begin{equation}\label{eq:planar-majorant-bound}
 0\le U^{(2)}_{nk}
 \le\frac{R^{-n}q^k}{1-\varepsilon R}.
\end{equation}
For $j\ge1$,
\[
 \binom{n-1}{j-1}\le\binom nj,\qquad
 0<\alpha_n^{(2)}\le\frac{1+\lambda}{1-\lambda},\qquad
 0<\beta_j^{(2)}\le\frac{1+\lambda\delta}{1-\lambda\delta}.
\]
For the rank-one term of \eqref{eq:K20},
$\varepsilon R<1$ and $q=T_\varepsilon(R)\ge\varepsilon$
give the separate bounds
\begin{equation}\label{eq:plane-rank-one-majorant}
 \varepsilon^n\le R^{-n},\qquad
 \varepsilon^k\le q^k,\qquad n,k\ge1.
\end{equation}
The estimates \eqref{eq:planar-majorant-bound} and
\eqref{eq:plane-rank-one-majorant} imply \eqref{eq:planar-hs}. Thus
\begin{equation}\label{eq:planar-HSsum}
 \sum_{n,k\ge1}|(K^{(2,m)}_{\lambda,\varepsilon})_{nk}|^2
 \le C^2\sum_{n,k\ge1}R^{-2n}q^{2k}<\infty.
\end{equation}
The majorant in \eqref{eq:planar-HSsum} is locally uniform in
$(\lambda,\varepsilon)$; dominated convergence proves continuity
in Hilbert--Schmidt norm. For $A=K^{(2,m)}_{\lambda,\varepsilon}A$,
Cauchy--Schwarz and \eqref{eq:planar-hs} give
\begin{equation}\label{eq:planar-fixed-decay}
 |A_n|\le CR^{-n}\|A\|_{\ell^2}
 \left(\sum_{k\ge1}q^{2k}\right)^{1/2},\qquad n\ge1.
\end{equation}
\end{proof}

The next proposition closes both directions of the PDE--operator correspondence and also covers the Neumann endpoint needed in the spectral-radius argument.
\begin{proposition}\label{prop:planar-equivalence}
Fix $0<\delta<1$, $0\le\varepsilon<1-\delta$ and
$0<\lambda<1$. The number $\lambda$ is a Steklov eigenvalue
of \eqref{eq:steklov} in the even (respectively odd) planar sector if and only if $1$ is an eigenvalue of $K^{(2,0)}_{\lambda,\varepsilon}$ (respectively $K^{(2,1)}_{\lambda,\varepsilon}$). The correspondence between the real eigenspace and the fixed-point kernel is a linear isomorphism. At $\lambda=0$, a nonzero fixed point of $K^{(2,0)}_{0,\varepsilon}=K^{(2,1)}_{0,\varepsilon}$ reconstructs a nonconstant harmonic Neumann function.
\end{proposition}
\begin{proof}
In the forward direction, smoothness, \Cref{lem:planar-split}
and \eqref{eq:planar-BC} yield a fixed point in $\ell^2$.
For each fixed $j$, the first series in \eqref{eq:planar-shifts}
is absolutely convergent by Cauchy--Schwarz: its weights
$\binom{k}{j}\varepsilon^{k-j}$ are square summable when
$0<\varepsilon<1$. At $\varepsilon=0$ the translation is the
identity; the second series is finite for each $n$.
If the extracted $A$ were zero, then $H_j=B_j=0$ and
\eqref{eq:planar-zeromode} would give $C=b=0$, so the original
function would vanish. Thus nonzero eigenfunctions give nonzero
fixed points. Conversely, let $A\in\ell^2$ satisfy the fixed-point equations
\eqref{eq:planar-fixed}. The bound \eqref{eq:planar-fixed-decay}
gives $|A_n|\le C_A R^{-n}$. For $1<S<R$ satisfying $S-\varepsilon>\delta$, the binomial identity in \eqref{eq:planar-shifts} gives
\begin{equation}\label{eq:planar-rebuild}
 |H_j|\le C\sum_{k\ge j}\bn{k}{j}\varepsilon^{k-j}S^{-k}
 =\frac{CS}{(S-\varepsilon)^{j+1}},\quad
 |B_j|\le C'\left(\frac{\delta^2}{S-\varepsilon}\right)^j.
\end{equation}
Set $\rho_0:=\delta^2/(S-\varepsilon)<\delta$.
For every compact $E\Subset\{|z|<S,\ |z-\varepsilon|>\rho_0\}$,
choose $r_+,r_-$ with
$0<\sup_E|z|\le r_+<S$ and
$\rho_0<r_-\le\inf_E|z-\varepsilon|$.
For every derivative order $N\ge0$, \eqref{eq:planar-rebuild}
gives the majorant
\begin{equation}\label{eq:planar-rebuild-normal}
 \begin{aligned}
 &\sum_{n\ge1}|A_n|\sup_E|\partial_z^Nz^n|
 +\sum_{j\ge1}|B_j|\sup_E
        |\partial_z^N(z-\varepsilon)^{-j}|\\
 &\quad\le C_{N,E}\sum_{n\ge1}(n+1)^N
                    \left(\frac{r_+}{S}\right)^n\\
 &\qquad{}+C_{N,E}\sum_{j\ge1}(j+N)^N
                  \left(\frac{\rho_0}{r_-}\right)^j<\infty.
 \end{aligned}
\end{equation}
Thus both series and their derivatives converge normally across
the two physical boundary circles. To justify changing the center
of the singular series there, \eqref{eq:planar-shifts} gives
\[
 \sum_{j=1}^n\binom{n-1}{j-1}
       \varepsilon^{n-j}|B_j|
 \le C\rho_0(\varepsilon+\rho_0)^{n-1},\qquad n\ge1.
\]
Since $\varepsilon+\rho_0<1$, the outer translated double
series is absolutely summable for radii
$r>\varepsilon+\rho_0$, including a neighborhood of $r=1$;
any fixed number of derivatives only adds polynomial factors.
Similarly the inner translated double series is bounded by
\[
 \sum_{k\ge1}|A_k|(\varepsilon+\rho)^k<\infty
 \qquad(\varepsilon+\rho<S),
\]
which holds on a neighborhood of $\rho=\delta$.
These estimates justify termwise translations in the boundary
residual identities.

For the even sector take $A_n\in\mathbb R$; for the odd sector
take $A_n\in i\mathbb R$. Define $B_j$ by the first (inner)
relation of \eqref{eq:planar-BC}, and $C,b$ by
\eqref{eq:planar-zeromode}
(with $C=b=0$ in the odd sector). Write $\mathcal B_{\mathrm i,j}$ and
$\mathcal B_{\mathrm o,n}$ for the positive Fourier coefficients,
as in \eqref{eq:planar-fourier-notation}, of
$2(\partial_\nu u-\lambda u)$ on the inner and outer circles.
For $j,n\ge1$, the Fourier residuals factor as
\begin{align}
 \mathcal B_{\mathrm i,j}
 &=-jH_j\delta^{j-1}+j\overline{B_j}\delta^{-j-1}
 -\lambda(H_j\delta^j+\overline{B_j}\delta^{-j})\\
 &=\frac{j-\lambda\delta}{\delta^{j+1}}
   \bigl(\overline{B_j}-\delta^{2j}\beta_j^{(2)}H_j\bigr)=0,
 \label{eq:planar-inner-residual}\\
 \mathcal B_{\mathrm o,n}
 &=(n-\lambda)A_n
 -(n+\lambda)\left(\overline{J_n}-\frac{b\varepsilon^n}{n}\right)\\
 &=(n-\lambda)\left[A_n-\alpha_n^{(2)}
 \left(\overline{J_n}-\frac{b\varepsilon^n}{n}\right)\right]=0,
 \label{eq:planar-outer-residual}
\end{align}
by \eqref{eq:planar-BC} and \eqref{eq:planar-fixed}.
Use the same Fourier normalization for degree zero:
$\mathcal B_{\mathrm o,0}$ and $\mathcal B_{\mathrm i,0}$ are the
means of $2(\partial_\nu u-\lambda u)$ on the respective circles.
Thus
\begin{equation}\label{eq:planar-zero-residual}
 \begin{aligned}
 \mathcal B_{\mathrm o,0}&=2(b-\lambda C)=0,\\
 \mathcal B_{\mathrm i,0}
 &=-\frac{2b}{\delta}
 -2\lambda\bigl(C+b\log\delta+\operatorname{Re}F(\varepsilon)\bigr)=0.
 \end{aligned}
\end{equation}
The negative Fourier coefficients are the conjugates of the
positive ones because the residual is real.
By \eqref{eq:planar-rebuild}, the boundary expansions and their
normal derivatives converge uniformly. The vanishing of every
Fourier coefficient in
\eqref{eq:planar-inner-residual}--\eqref{eq:planar-zero-residual}
therefore yields
\begin{equation}\label{eq:planar-reconstruction-boundary}
 (\partial_\nu u-\lambda u)|_{\partial\Omega_{
\delta,\varepsilon}^{(2)}}=0.
\end{equation}
If $u\equiv0$, the uniqueness part of \Cref{lem:planar-split}
forces $A_n=0$ for all $n\ge1$. Starting instead from an
eigenfunction, the relations \eqref{eq:planar-BC} and
\eqref{eq:planar-zeromode} determine its $B_j,b,C$ uniquely
from $A$. Reconstruction thus returns that eigenfunction.
For the odd sector the real sequence used by $K^{(2,1)}$ is
$A^{\mathrm I}$ and the complex regular coefficients are
$A=iA^{\mathrm I}$; the even sector uses $A=A^{\mathrm R}$.
Both identifications are real-linear. Green's identity, as in
the higher-dimensional proof, gives exactly the weak formulation
\eqref{eq:eigenspace-definition}.

The correspondence is also bounded in the same norms as above.
Let $u$ be a real harmonic $H^1$ function in the planar shell. On a fixed auxiliary inner circle, interior
estimates for $g=2\partial_z u$ bound $M_{\mathrm i}$ by
$C\|u\|_{H^1}$. The Cauchy bounds
\eqref{eq:planar-h-est} and \eqref{eq:planar-coeff-deriv} then
bound $b\log|z-\varepsilon|+\operatorname{Re}G(z)$ in every
smooth norm near the outer circle, since
$r_{\mathrm i}<1-\varepsilon$. Subtracting this function from
$u$ gives the $H^{1/2}$ outer trace of
$C+\operatorname{Re}F$. Its nonzero Fourier coefficients are
$A_n/2$ and $\overline{A_n}/2$, so
\[
 \sum_{n\ge1}(1+n^2)^{1/2}|A_n|^2
 \le C_{\delta,\varepsilon}\|u\|_{H^1}^2.
\]
Thus coefficient extraction into $\ell^2$ is bounded, including
on each parity sector. The planar estimates
\eqref{eq:planar-fixed-decay} and \eqref{eq:planar-rebuild-normal},
with the zero coefficients from \eqref{eq:planar-zeromode}, give
$\|u\|_{C^p}\le C_p\|A\|_{\mathscr X_R}$ for arbitrary
$A\in\mathscr X_R$, and
$\|u\|_{H^1}\le C\|A\|_{\ell^2}$ on the fixed-point kernel.
The claimed real-linear isomorphism is therefore topological.
At $\lambda=0$, the zero-mode equations require $b=0$ and
leave $C$ free. Taking $C=b=0$ and using
\eqref{eq:planar-rebuild} gives a nonconstant harmonic function
with zero normal derivative whenever $A\ne0$ is a fixed point.
Nonconstancy follows from uniqueness in \Cref{lem:planar-split}:
a constant has all $A_n=B_j=0$ and $b=0$.
The continuous extension of the positive-parameter formula
$C=-\operatorname{Re}F(\varepsilon)/D_\lambda^{(2)}$ merely
chooses a different additive constant at $\lambda=0$.
At $\lambda=0$, the constant coefficient is free, whereas the
fixed-point kernel describes the nonconstant modes.
\end{proof}

\subsection{Spectral-radius crossing}\label{subsec:plane-spectral}
By \eqref{eq:K21}--\eqref{eq:K20}, the difference
$K^{(2,0)}-K^{(2,1)}$ is an entrywise positive rank-one operator
when $\lambda,\varepsilon>0$. We use
\Cref{lem:prelim-perron,lem:prelim-radius-order} to compare the
two spectral radii.
\begin{proposition}\label{prop:planar-crossing}
Fix $0<\delta<1$, $0\le\varepsilon<1-\delta$, and set
$\mathcal R_2(\lambda,\varepsilon)
=\spr(K^{(2,0)}_{\lambda,\varepsilon})$ for
$0\le\lambda<1$. The function $\lambda\mapsto\mathcal R_2(\lambda,\varepsilon)$ is continuous and strictly increasing on $[0,1)$, with
\begin{equation}\label{eq:planar-cross}
 \mathcal R_2(0,\varepsilon)<1,
 \qquad\lim_{\lambda\uparrow1}\mathcal R_2(\lambda,\varepsilon)=+\infty.
\end{equation}
Its unique level-one solution equals $\sigma_1^{(2)}(\delta,\varepsilon)$. For fixed $0<\lambda<1$, $\mathcal R_2(\lambda,\varepsilon)$ is strictly increasing in $\varepsilon$.
\end{proposition}
\begin{proof}
For $\varepsilon>0$, the $j=1$ summand in \eqref{eq:K21} gives
\[
 (K^{(2,0)})_{nk}\ge(K^{(2,1)})_{nk}
 \ge\alpha_n^{(2)}k\delta^2\beta_1^{(2)}
       \varepsilon^{n+k-2}>0.
\]
Together with compactness from \Cref{lem:planar-compact}, this
allows us to apply \Cref{lem:prelim-perron}.
The multiplier identities in \eqref{eq:planar-factors} give
\begin{equation}\label{eq:planar-parameter-derivatives}
 \begin{aligned}
 \partial_\lambda\alpha_n^{(2)}
 &=\frac{2n}{(n-\lambda)^2}>0,\\
 \partial_\lambda\beta_j^{(2)}
 &=\frac{2j\delta}{(j-\lambda\delta)^2}>0,\\
 \partial_\lambda\left(\frac{\lambda}{D_\lambda^{(2)}}\right)
 &=\frac{1+\delta^{-1}}{(D_\lambda^{(2)})^2}>0.
 \end{aligned}
\end{equation}
For $0\le\lambda_1<\lambda_2<1$ let $v_1>0$ and $w_2>0$ be,
respectively, a right Perron vector of the smaller operator and a
left Perron vector of the larger one. With
$r_i=\mathcal R_2(\lambda_i,\varepsilon)$,
\begin{equation}\label{eq:planar-lambda-order}
 (r_2-r_1)\langle w_2,v_1\rangle
 =\left\langle w_2,
 (K^{(2,0)}_{\lambda_2,\varepsilon}
 -K^{(2,0)}_{\lambda_1,\varepsilon})v_1\right\rangle>0.
\end{equation}
Thus $r_2>r_1$.
At $\varepsilon=0$ the matrices are diagonal with entries
$\alpha_n^{(2)}\beta_n^{(2)}\delta^{2n}$. For fixed
$0\le\lambda<1$, both factors
\[
 \alpha_n^{(2)}=1+\frac{2\lambda}{n-\lambda},\qquad
 \beta_n^{(2)}=1+\frac{2\lambda\delta}{n-\lambda\delta}
\]
are nonincreasing in $n$, and $\delta^{2n}$ strictly decreases.
The unique maximizing coordinate is therefore $n=1$.
Every diagonal entry increases strictly with $\lambda$;
therefore so does the spectral radius.
The spectral radius remains positive because
$\mathcal R_2\ge(K^{(2,0)})_{11}>0$.
Joint continuity in $(\lambda,\varepsilon)$ follows from
\eqref{eq:compact-radius-continuity} and the Hilbert--Schmidt-norm
continuity in \Cref{lem:planar-compact}.

At $\varepsilon=\lambda=0$, the diagonal operator in
\eqref{eq:K21} gives
$\mathcal R_2(0,0)=\delta^2<1$. If
$\mathcal R_2(0,\varepsilon_*)=1$ for an admissible
$\varepsilon_*>0$, the Perron fixed point and
\Cref{prop:planar-equivalence} would produce a nonconstant
harmonic function $u$ satisfying $\partial_\nu u=0$. By
Green's identity,
\begin{equation}\label{eq:plane-neumann}
 \int_{\Omega^{(2)}_{\delta,\varepsilon_*}}|\nabla u|^2
 =\int_{\partial\Omega^{(2)}_{\delta,\varepsilon_*}}u\partial_\nu u=0,
\end{equation}
a contradiction. If the radius were greater than or equal to
one at some admissible displacement, continuity on the segment
from zero to that displacement, starting at $\delta^2<1$,
would produce just such a level-one value. Hence the first
assertion in \eqref{eq:planar-cross} holds. The diagonal spectral-radius estimate
\eqref{eq:diagonal-spr-bound} and the $n=k=j=1$ summand give
\begin{equation}\label{eq:planar-diverge}
 \mathcal R_2(\lambda,\varepsilon)
 \ge (K^{(2,0)}_{\lambda,\varepsilon})_{11}
 \ge \alpha_1^{(2)}\beta_1^{(2)}\delta^2\longrightarrow\infty
 \quad(\lambda\uparrow1).
\end{equation}
Together with the strict increase in $\lambda$ proved in \eqref{eq:planar-lambda-order}, \eqref{eq:planar-cross} yields a unique root $\lambda_*\in(0,1)$ of $\mathcal R_2(\lambda_*,\varepsilon)=1$.

For $\varepsilon>0$ and $\lambda>0$, the rank-one term
in \eqref{eq:K20} has strictly positive entries. Let $v_1>0$ be
a right Perron vector of $K^{(2,1)}$ and $w_0>0$ a left Perron vector
of $K^{(2,0)}$. With $r_m=\spr(K^{(2,m)}_{\lambda,\varepsilon})$,
\begin{equation}\label{eq:planar-branchpair}
 (r_0-r_1)\langle w_0,v_1\rangle
 =\langle w_0,
 (K^{(2,0)}_{\lambda,\varepsilon}-
 K^{(2,1)}_{\lambda,\varepsilon})v_1\rangle>0.
\end{equation}
The right-hand side of \eqref{eq:planar-branchpair} equals
\[
 \frac{\lambda}{D_\lambda^{(2)}}
 \left(\sum_{n\ge1}(w_0)_n\frac{\alpha_n^{(2)}}{n}
                    \varepsilon^n\right)
 \left(\sum_{k\ge1}(v_1)_k\varepsilon^k\right)>0.
\]
Both series converge absolutely by Cauchy--Schwarz, because
$w_0,v_1\in\ell^2$ and $0<\varepsilon<1$.
The denominator $\langle w_0,v_1\rangle$ is finite and positive.
Consequently,
\begin{equation}\label{eq:planar-branchorder}
 \spr(K^{(2,1)}_{\lambda,\varepsilon})
 <\spr(K^{(2,0)}_{\lambda,\varepsilon}).
\end{equation}
For $\varepsilon=0$, \eqref{eq:K20} gives
$K^{(2,0)}_{\lambda,0}=K^{(2,1)}_{\lambda,0}$.
By \Cref{prop:planar-equivalence}, the level-one equation and
\eqref{eq:planar-branchorder} yield
\begin{equation}\label{eq:planar-root-exclusion}
 \begin{gathered}
 \spr(K^{(2,0)}_{\lambda_*,\varepsilon})=1,\\
 0<\lambda<\lambda_*\ \Longrightarrow\quad
 \max_{m=0,1}\spr(K^{(2,m)}_{\lambda,\varepsilon})<1
 \ \Longrightarrow\ \mathscr E_\lambda^{(2)}(\delta,\varepsilon)=\{0\},\\
 \lambda_*=\sigma_1^{(2)}(\delta,\varepsilon).
 \end{gathered}
\end{equation}
For $0\le\varepsilon_1<\varepsilon_2<1-\delta$ and
$0<\lambda<1$, \eqref{eq:K20} gives
\begin{equation}\label{eq:planar-ecc-entry-order}
 K^{(2,0)}_{\lambda,\varepsilon_2}
 >K^{(2,0)}_{\lambda,\varepsilon_1}
 \quad\text{entrywise}.
\end{equation}
Indeed, the $(n,k)$ entry of this difference is bounded below by
\[
 \frac{\alpha_n^{(2)}\lambda}{nD_\lambda^{(2)}}
 (\varepsilon_2^{n+k}-\varepsilon_1^{n+k})>0,
\]
including $\varepsilon_1=0$. Take a nonnegative right eigenvector
$v_1$ of the smaller operator at its spectral radius
(specifically $v_1=e_1$ if $\varepsilon_1=0$) and a strictly
positive left Perron vector $w_2$ of the larger operator.
The absolute convergence and strict positivity of the pairing
follow from \Cref{lem:prelim-radius-order}.
Writing $r_i=\mathcal R_2(\lambda,\varepsilon_i)$, the equality
\begin{equation}\label{eq:planar-ecc-pair}
 (r_2-r_1)\langle w_2,v_1\rangle
 =\langle w_2,(K^{(2,0)}_{\lambda,\varepsilon_2}
 -K^{(2,0)}_{\lambda,\varepsilon_1})v_1\rangle>0
\end{equation}
proves strict eccentricity monotonicity of $\mathcal R_2$.
\end{proof}

\begin{proof}[Proof of \Cref{thm:main} when $d=2$]
Set $\lambda_i=\sigma_1^{(2)}(\delta,\varepsilon_i)$.
The strict eccentricity comparison and the unique crossing in
\Cref{prop:planar-crossing} give
\begin{equation}\label{eq:planar-main-final}
 \mathcal R_2(\lambda_1,\varepsilon_1)=1,
 \qquad \mathcal R_2(\lambda_1,\varepsilon_2)>1,
 \qquad \mathcal R_2(\lambda_2,\varepsilon_2)=1
 \quad\Longrightarrow\quad \lambda_2<\lambda_1.
\end{equation}
For $\varepsilon>0$, the even Perron eigenvalue is algebraically
simple, and \eqref{eq:planar-branchorder} excludes the odd branch
at the root. Hence the first positive eigenspace is one-dimensional.

At $\varepsilon=0$ both parity operators are diagonal, with entries
\begin{equation}\label{eq:plane-concentric}
 \eta_n^{(2)}(\lambda)=\alpha_n^{(2)}\beta_n^{(2)}\delta^{2n}
 =\frac{(n+\lambda)(n+\lambda\delta)}
 {(n-\lambda)(n-\lambda\delta)}\delta^{2n},\qquad n\ge1.
\end{equation}
For fixed $0<\lambda<1$,
\begin{equation}\label{eq:plane-diag-order}
 \frac{n+1+\lambda}{n+1-\lambda}
 <\frac{n+\lambda}{n-\lambda},\qquad
 \frac{n+1+\lambda\delta}{n+1-\lambda\delta}
 <\frac{n+\lambda\delta}{n-\lambda\delta},
 \qquad 0<\delta^2<1.
\end{equation}
Consequently $\eta_{n+1}^{(2)}(\lambda)<\eta_n^{(2)}(\lambda)$,
so among nonconstant angular modes the first level-one crossing
occurs at $n=1$. For $\varepsilon=0$ and a nonconstant radial eigenfunction
$u(r)=C+b\log r$ with $b\ne0$, the degree-zero boundary equations
\eqref{eq:planar-zeromode} become
\begin{equation}\label{eq:planar-radial-system}
 b=\lambda C,\qquad
 -\frac b\delta=\lambda(C+b\log\delta)
 \quad\Longrightarrow\quad
 \lambda\log(\delta^{-1})=1+\delta^{-1}.
\end{equation}
Hence the radial $n=0$ positive eigenvalue is
\begin{equation}\label{eq:planar-radial-positive}
 \lambda_{0,+}^{(2)}
 =\frac{1+\delta^{-1}}{\log(\delta^{-1})}>1,
\end{equation}
where the last inequality follows from $\log t<t-1$ for $t>1$.
Since \eqref{eq:range} gives $\sigma_1^{(2)}<1$, the radial
positive mode cannot be first. At $n=1$ the angular eigenspace is
\begin{equation}\label{eq:planar-first-angular-span}
 \operatorname{span}_{\RR}\{\cos\theta,\sin\theta\},
 \qquad \dim_{\RR}=2.
\end{equation}
The two radial boundary conditions determine a one-dimensional
radial eigenspace for each angular function, proving multiplicity two.
\end{proof}

\section{Spectral characterization and concentric perturbations}\label{sec:all-dimensional}
\subsection{Spectral-radius characterization}
For $d\ge3$, let $K^{(d,0)}_{\lambda,\varepsilon}$ be defined by
\eqref{eq:Kd0}, and for $d=2$ by \eqref{eq:K20}.

\begin{corollary}\label{cor:all-dim-crossing}
Fix $d\ge2$ and $0<\delta<1$. The function
\[
 \mathcal R_d(\lambda,\varepsilon)
 :=\spr(K^{(d,0)}_{\lambda,\varepsilon}),
 \qquad 0\le\lambda<1,\quad 0\le\varepsilon<1-\delta,
\]
is jointly continuous, strictly increasing in $\lambda$, and has a
unique level-one root. Moreover,
\begin{equation}\label{eq:intro-crossing}
 \mathcal R_d\bigl(\sigma_1^{(d)}(\delta,\varepsilon),\varepsilon\bigr)=1,
 \qquad
 \mathcal R_d(\lambda,\varepsilon_2)>
 \mathcal R_d(\lambda,\varepsilon_1)
\end{equation}
whenever $0<\lambda<1$ and
$0\le\varepsilon_1<\varepsilon_2<1-\delta$.
The map $\varepsilon\mapsto\sigma_1^{(d)}(\delta,\varepsilon)$
is continuous on $[0,1-\delta)$.
\end{corollary}
\begin{proof}
The estimates in \Cref{prop:multi-crossing,prop:planar-crossing} give,
for every admissible $\varepsilon$,
\begin{equation}\label{eq:all-dim-crossing-summary}
\begin{gathered}
 \mathcal R_d(0,\varepsilon)<1,
 \qquad \lim_{\lambda\uparrow1}\mathcal R_d(\lambda,\varepsilon)=+\infty,\\
 0\le\lambda_1<\lambda_2<1
 \quad\Longrightarrow\quad
 \mathcal R_d(\lambda_1,\varepsilon)
 <\mathcal R_d(\lambda_2,\varepsilon).
\end{gathered}
\end{equation}
The fixed-point equivalences
\Cref{prop:multi-equivalence,prop:planar-equivalence} and the angular
comparisons in \Cref{prop:multi-identify,prop:planar-crossing} yield
\[
 \exists!\,\lambda_*(\varepsilon)\in(0,1):\quad
 \mathcal R_d(\lambda_*(\varepsilon),\varepsilon)=1,
 \qquad
 \lambda_*(\varepsilon)=\sigma_1^{(d)}(\delta,\varepsilon).
\]
The entrywise comparisons
\eqref{eq:multi-ecc-pairing} and \eqref{eq:planar-ecc-pair} give
\[
 \bigl[\mathcal R_d(\lambda,\varepsilon_2)
  -\mathcal R_d(\lambda,\varepsilon_1)\bigr]
 \langle w_2,v_1\rangle>0,
 \qquad 0<\lambda<1,\quad\varepsilon_2>\varepsilon_1,
\]
with $v_1\ge0$ nonzero and $w_2>0$. In particular,
for $0\le\varepsilon_1<\varepsilon_2<1-\delta$,
\[
 1=\mathcal R_d(\lambda_*(\varepsilon_1),\varepsilon_1)
 <\mathcal R_d(\lambda_*(\varepsilon_1),\varepsilon_2),
 \qquad
 \mathcal R_d(\lambda_*(\varepsilon_2),\varepsilon_2)=1,
\]
so strict increase in $\lambda$ yields
$\lambda_*(\varepsilon_2)<\lambda_*(\varepsilon_1)$.

Fix $\varepsilon_0\in[0,1-\delta)$ and choose
$\lambda_-<\lambda_*(\varepsilon_0)<\lambda_+$ with
$0<\lambda_-<\lambda_+<1$. By
\eqref{eq:all-dim-crossing-summary},
\[
 \mathcal R_d(\lambda_-,\varepsilon_0)<1
 <\mathcal R_d(\lambda_+,\varepsilon_0).
\]
Joint continuity preserves these inequalities for
$\varepsilon$ near $\varepsilon_0$; uniqueness of the root gives
\[
 \lambda_-<\sigma_1^{(d)}(\delta,\varepsilon)<\lambda_+,
 \qquad |\varepsilon-\varepsilon_0|<h
\]
for sufficiently small $h>0$. Letting
$\lambda_-,\lambda_+\to\sigma_1^{(d)}(\delta,\varepsilon_0)$,
we obtain
\begin{equation}\label{eq:all-dim-root-continuity}
 \lim_{\substack{\varepsilon\to\varepsilon_0\\
                    0\le\varepsilon<1-\delta}}
 \sigma_1^{(d)}(\delta,\varepsilon)
 =\sigma_1^{(d)}(\delta,\varepsilon_0).
\end{equation}
\end{proof}

\subsection{The concentric shell}
\begin{corollary}\label{cor:concentric-first-root}
For every integer $d\ge2$ and $0<\delta<1$,
$\sigma_1^{(d)}(\delta,0)$ is the unique $\lambda\in(0,1)$
satisfying
\begin{equation}\label{eq:all-dim-concentric-scalar}
 \delta^d\frac{(d-1+\lambda)(1+\lambda\delta)}
 {(1-\lambda)(d-1-\lambda\delta)}=1.
\end{equation}
The corresponding eigenspace has dimension $d$.
\end{corollary}
\begin{proof}
By \eqref{eq:concentric-highdim-mode-exclusion} and
\eqref{eq:radial-zero-d} for $d\ge3$, and by
\eqref{eq:plane-diag-order} and \eqref{eq:planar-radial-positive}
for $d=2$, the first positive mode has spherical-harmonic degree one.
For a nonzero degree-one spherical harmonic
$Y(\widehat x)=b\cdot\widehat x$ with $b\in\mathbb R^d$, write
\begin{equation}\label{eq:concentric-direct-degree-one}
 u(r,\widehat x)=(A_1r+B_1r^{-(d-1)})Y(\widehat x),
 \qquad \widehat x\in\mathbb S^{d-1}.
\end{equation}
The outer and inner boundary conditions give
\begin{equation}\label{eq:concentric-direct-radial-system}
 \begin{pmatrix}1-\lambda&-(d-1+\lambda)\\
 -\delta^d(1+\lambda\delta)&d-1-\lambda\delta
 \end{pmatrix}
 \begin{pmatrix}A_1\\ B_1\end{pmatrix}=0.
\end{equation}
A nontrivial solution exists precisely when the determinant vanishes,
which is \eqref{eq:all-dim-concentric-scalar}.
Let $L(\lambda)$ denote the left-hand side of this equation. Then
\begin{equation}\label{eq:all-dim-concentric-derivative}
 \begin{gathered}
 L(0)=\delta^d<1,\qquad
 \lim_{\lambda\uparrow1}L(\lambda)=+\infty,\\
 \frac{L'(\lambda)}{L(\lambda)}
 =\frac1{d-1+\lambda}+\frac1{1-\lambda}
 +\frac{\delta}{1+\lambda\delta}
 +\frac{\delta}{d-1-\lambda\delta}>0.
 \end{gathered}
\end{equation}
Thus $L(\lambda)=1$ has exactly one root in $(0,1)$.
For each $Y$, \eqref{eq:concentric-direct-radial-system} has a
one-dimensional kernel. The degree-one spherical harmonics form
\[
 \dim_{\mathbb R}
 \operatorname{span}\{\widehat x_1,\ldots,\widehat x_d\}=d,
\]
so the eigenspace has dimension $d$.
\end{proof}

\subsection{Quadratic splitting near concentricity}
We use the diagonal quantities $\eta_n^{(d)}$ defined in
\eqref{eq:concentric-highdim-eta} for $d\ge3$ and in
\eqref{eq:plane-concentric} for $d=2$.

\begin{proposition}\label{prop:local-concentric-splitting}
Fix $d\ge2$, $0<\delta<1$, and put
$\lambda_*:=\sigma_1^{(d)}(\delta,0)$.
Let $\lambda_{\mathrm a}(\varepsilon)$ denote the lowest positive
axial eigenvalue (even sector for $d=2$), and let
$\lambda_{\mathrm t}(\varepsilon)$ be the lowest eigenvalue
in a fixed first transverse component (odd sector for $d=2$).
For sufficiently small $\varepsilon\ge0$, these branches have even
real-analytic extensions to signed $\varepsilon$ near zero. Define
\begin{equation}\label{eq:local-concentric-constants}
 \begin{aligned}
 g_d&:=\frac{\alpha_2^{(d)}(\lambda_*)}
 {\alpha_1^{(d)}(\lambda_*)[1-\eta_2^{(d)}(\lambda_*)]}>0,\\
 p_d&:=\begin{cases}
 \displaystyle\frac{\alpha_1^{(2)}(\lambda_*)\lambda_*}
 {D_{\lambda_*}^{(2)}},&d=2,\\[5pt]
 \displaystyle\frac{(d-2)\alpha_1^{(d)}(\lambda_*)
            \lambda_*\delta^{d-1}}{D_{\lambda_*}^{(d)}},&d\ge3.
 \end{cases}
 \end{aligned}
\end{equation}
Then $p_d>0$ and, with the derivative in $\lambda$,
\begin{equation}\label{eq:local-concentric-expansion}
 \begin{aligned}
 \lambda_{\mathrm a}(\varepsilon)
 &=\lambda_*
 -\frac{2(d-1)g_d+p_d}{(\eta_1^{(d)})'(\lambda_*)}
       \varepsilon^2+O(\varepsilon^4),\\
 \lambda_{\mathrm t}(\varepsilon)
 &=\lambda_*
 -\frac{d g_d}{(\eta_1^{(d)})'(\lambda_*)}
       \varepsilon^2+O(\varepsilon^4),\\
 \lambda_{\mathrm t}(\varepsilon)-\lambda_{\mathrm a}(\varepsilon)
 &=\frac{(d-2)g_d+p_d}{(\eta_1^{(d)})'(\lambda_*)}
       \varepsilon^2+O(\varepsilon^4).
 \end{aligned}
\end{equation}
The denominator is strictly positive, and the $O(\varepsilon^4)$
bounds are uniform on a sufficiently small interval, with constants
depending on $d,\delta$. In particular,
\[
 \lambda_{\mathrm a}'(0)=\lambda_{\mathrm t}'(0)=0,
 \qquad \lambda_{\mathrm a}''(0)<0,\quad
 \lambda_{\mathrm t}''(0)<0.
\]
The axial branch satisfies
$\lambda_{\mathrm a}(\varepsilon)=\sigma_1^{(d)}(\delta,\varepsilon)$.
For $d\ge3$, the transverse branch equals
$\tau_1^{(d)}(\delta,\varepsilon)$; for $d=2$ it is the lowest
odd-sector eigenvalue.
\end{proposition}
\begin{proof}
For $m=0,1$, abbreviate
$K_m(\lambda,\varepsilon)=K^{(d,m)}_{\lambda,\varepsilon}$,
and complexify $\ell^2$. Choose a complex neighborhood of
$\lambda_*$ avoiding the poles of the Robin multipliers and
$1/D_\lambda^{(d)}$, and $0<t<1-\delta$.
For $|\varepsilon|\le t$, choose $R>1$ with
$tR<1$ and $q=T_t(R)<1$. The majorants in
\Cref{lem:multi-compact,lem:planar-compact} yield
\begin{equation}\label{eq:local-hs-holomorphy}
 \begin{aligned}
 |(K_m(\lambda,\varepsilon))_{nk}|
 &\le C(k+1)^{d-2}R^{-n}q^k,\\
 \sum_{n,k\ge1}|(K_m(\lambda,\varepsilon))_{nk}|^2
 &\le C^2\left(\sum_{n\ge1}R^{-2n}\right)
 \left(\sum_{k\ge1}(k+1)^{2(d-2)}q^{2k}\right)<\infty.
 \end{aligned}
\end{equation}
The finite coordinate compressions converge uniformly in
Hilbert--Schmidt norm. As their entries are holomorphic, so is
$(\lambda,\varepsilon)\mapsto K_m(\lambda,\varepsilon)$ in
operator norm.

Let $e_1=(1,0,\ldots)$ and $Q=\Id-e_1\otimes e_1$.
At $(\lambda_*,0)$,
\begin{equation}\label{eq:local-complement-gap}
 \begin{gathered}
 K_m(\lambda_*,0)
 =\operatorname{diag}(1,\eta_2^{(d)}(\lambda_*),\ldots),\\
 \|QK_m(\lambda_*,0)Q\|=\eta_2^{(d)}(\lambda_*)<1,\\
 (\Id-QK_mQ)^{-1}=\sum_{j=0}^\infty(QK_mQ)^j
 \quad\text{near }(\lambda_*,0).
 \end{gathered}
\end{equation}
Set $A=e_1+y_m$ with $QA=y_m$. The fixed-point equation
$K_mA=A$ is equivalent to
\begin{equation}\label{eq:local-schur-complement}
 \begin{aligned}
 y_m(\lambda,\varepsilon)
 &= (\Id-QK_mQ)^{-1}QK_me_1,\\
 f_m(\lambda,\varepsilon)
 &:=\langle e_1,K_m(e_1+y_m)\rangle-1=0.
 \end{aligned}
\end{equation}
Here $A_1\ne0$ for every nonzero fixed point because
$\|QK_mQ\|<1$. At $\varepsilon=0$,
$f_m(\lambda,0)=\eta_1^{(d)}(\lambda)-1$, so
\begin{equation}\label{eq:local-root-derivative}
 \partial_\lambda f_m(\lambda_*,0)
 =(\eta_1^{(d)})'(\lambda_*)
 =\frac1{d-1+\lambda_*}+\frac1{1-\lambda_*}
 +\frac{\delta}{1+\lambda_*\delta}
 +\frac{\delta}{d-1-\lambda_*\delta}>0.
\end{equation}
The analytic implicit-function theorem gives local real-analytic
roots $\lambda_m(\varepsilon)$ with $\lambda_m(0)=\lambda_*$. Let
$V=\operatorname{diag}((-1)^n)_{n\ge1}$. The coefficient parity yields
\begin{equation}\label{eq:local-parity}
 V^2=\Id,\quad Ve_1=-e_1,\quad
 VK_m(\lambda,\varepsilon)V=K_m(\lambda,-\varepsilon),
 \quad f_m(\lambda,-\varepsilon)=f_m(\lambda,\varepsilon).
\end{equation}
Uniqueness in the implicit-function theorem now gives
$\lambda_m(-\varepsilon)=\lambda_m(\varepsilon)$.

In operator norm, expand at $(\lambda_*,0)$:
\[
 K_m(\lambda_*,\varepsilon)
 =K_m(\lambda_*,0)+\varepsilon L_m+\varepsilon^2N_m
 +O(\varepsilon^3).
\]
Since the entry $(n,k)$ has powers $\varepsilon^{n+k-2j}$,
the first derivative has nonzero $(1,k)$ and $(k,1)$ entries
only when $k=2$. More precisely,
\begin{equation}\label{eq:local-adjacent-indices}
 QL_me_1=(L_m)_{21}e_2,
 \qquad \langle e_1,L_mQz\rangle=(L_m)_{12}\langle e_2,z\rangle,
 \quad z\in\ell^2.
\end{equation}
Direct extraction from \eqref{eq:Kd0}--\eqref{eq:Kd1} and
\eqref{eq:K21}--\eqref{eq:K20} gives
\begin{equation}\label{eq:local-matrix-coefficients}
 \begin{array}{c|ccc}
 & (L_m)_{12}&(L_m)_{21}&(N_m)_{11}\\ \hline
 d\ge3,\ m=0
 &(d-1)&2\alpha_2^{(d)}/\alpha_1^{(d)}&p_d\\
 d\ge3,\ m=1
 &d&\alpha_2^{(d)}/\alpha_1^{(d)}&0\\
 d=2,\ m=0
 &2&\alpha_2^{(2)}/\alpha_1^{(2)}&p_2\\
 d=2,\ m=1
 &2&\alpha_2^{(2)}/\alpha_1^{(2)}&0
 \end{array}
\end{equation}
where all multipliers are evaluated at $\lambda_*$.
The quadratic coefficient of the Schur complement is
\begin{equation}\label{eq:local-schur-quadratic}
 \begin{aligned}
 [\varepsilon^2]f_m(\lambda_*,\varepsilon)
 &=(N_m)_{11}
 +\big\langle e_1,L_mQ(\Id-QK_m(\lambda_*,0)Q)^{-1}
 QL_me_1\big\rangle\\
 &=(N_m)_{11}
 +\frac{(L_m)_{12}(L_m)_{21}}
 {1-\eta_2^{(d)}(\lambda_*)}
 =:c_m
 =\begin{cases}
 2(d-1)g_d+p_d,&m=0,\\
 dg_d,&m=1.
 \end{cases}
 \end{aligned}
\end{equation}
By joint analyticity and \eqref{eq:local-parity},
\begin{equation}\label{eq:local-scalar-taylor}
 \begin{aligned}
 f_m(\lambda,\varepsilon)
 &=(\eta_1^{(d)})'(\lambda_*)(\lambda-\lambda_*)
   +c_m\varepsilon^2\\
 &\quad+O\!\left(|\lambda-\lambda_*|^2
 +|\lambda-\lambda_*|\,|\varepsilon|^2
 +|\varepsilon|^4\right),\\
 \lambda_m(\varepsilon)-\lambda_*
 &=-\frac{c_m}{(\eta_1^{(d)})'(\lambda_*)}
 \varepsilon^2+O(\varepsilon^4).
 \end{aligned}
\end{equation}
Since $g_d>0$ and $p_d>0$,
\[
 c_0-c_1=(d-2)g_d+p_d>0,
\]
which yields \eqref{eq:local-concentric-expansion}.
Even analyticity gives the stated $O(\varepsilon^4)$ bound
on a sufficiently small interval.

For real $\varepsilon\ge0$ at a local root, the Neumann series in
\eqref{eq:local-complement-gap} is nonnegative entrywise; hence
$A_m=e_1+y_m\ge0$ is nonzero. For $\varepsilon>0$, let
$w_m>0$ be a left Perron vector. Then
\begin{equation}\label{eq:local-root-radius}
 K_m(\lambda_m(\varepsilon),\varepsilon)A_m=A_m,
 \quad
 [\spr(K_m(\lambda_m(\varepsilon),\varepsilon))-1]
 \langle w_m,A_m\rangle=0,
 \quad \spr(K_m(\lambda_m(\varepsilon),\varepsilon))=1.
\end{equation}
At $\varepsilon=0$ the last equality follows from the diagonal
matrix. The parameter monotonicity established in
\Cref{prop:multi-crossing,prop:planar-crossing}, and for the first
transverse sector by the same strictly increasing Robin multipliers,
gives
\begin{equation}\label{eq:local-lowest-sector}
 0<\lambda<\lambda_m(\varepsilon)
 \quad\Longrightarrow\quad
 \spr(K_m(\lambda,\varepsilon))<1
 \quad\Longrightarrow\quad
 \ker(\Id-K_m(\lambda,\varepsilon))=\{0\}.
\end{equation}
The reconstruction equivalences in
\Cref{prop:multi-equivalence,prop:planar-equivalence} identify these
roots as the lowest positive eigenvalues in the two sectors.
For sufficiently small $\varepsilon>0$,
\[
 \begin{cases}
 \lambda_{\mathrm a}(\varepsilon)
 <\lambda_{\mathrm t}(\varepsilon)=\tau_1^{(d)}(\delta,\varepsilon)
 <\tau_m^{(d)}(\delta,\varepsilon),& d\ge3,\ m\ge2,\\[2pt]
 \lambda_{\mathrm a}(\varepsilon)<\lambda_{\mathrm t}(\varepsilon),
 &d=2.
 \end{cases}
\]
The angular-sector decomposition and \Cref{lem:higher-angular}
therefore give
$\lambda_{\mathrm a}(\varepsilon)
 =\sigma_1^{(d)}(\delta,\varepsilon)$.
At $\varepsilon=0$, this follows from \Cref{cor:concentric-first-root}.
\end{proof}

\section{Concluding remarks}\label{sec:conclusions}
The proof of strict eccentricity monotonicity combines the
spectral-radius comparison of positive compact coefficient operators
with the angular-sector ordering. The rank-one terms in
\eqref{eq:Kd0} and \eqref{eq:K20} come from eliminating the
degree-zero multipoles, $r^{-(d-2)}$ for $d\ge3$ and $\log r$
for $d=2$. The local expansion in
\Cref{prop:local-concentric-splitting} further quantifies the
quadratic decrease of the first eigenvalue and the splitting of
its first eigenspace. The contact limit and extensions beyond
spherical boundaries are natural directions for further study.

\paragraph{Declaration of AI assistance.}

The authors used ChatGPT (OpenAI) for exploratory mathematical discussions and editorial assistance during the preparation of this manuscript. The authors take full responsibility for all mathematical statements, arguments, and proofs.

\end{document}